\documentclass[12pt]{article}

\usepackage{color}
\usepackage{amsmath}
\usepackage{amssymb}
\usepackage{amsbsy}

\usepackage{amsfonts}
\usepackage{tikz}

\usepackage[margin=1in]{geometry}

\newtheorem{lem}{Lemma}[section]%
\newtheorem{thm}[lem]{Theorem}%
\newtheorem{defi}[lem]{Definition}%
\newtheorem{prob}{Problem}%
\newtheorem{prop}[lem]{Proposition}%
\newtheorem{rem}[lem]{Remark}%

\def\a{\alpha}    
 \def\s{\sigma}

\def\G{\Gamma}

   \def\olF{\overline F}
 \def\lg{\langle} \def\rg{\rangle}
\def\nd{\mathrel{\bigm|\kern-.7em/}}

\def\f{\noindent}

 \def\GL{\hbox{\rm GL}}

\DeclareMathOperator{\Aut}{Aut}

\def\soc{\hbox{\rm soc}}

\def\Cay{\hbox{\rm Cay}}

\def\mod{\hbox{\rm mod }}

\def\K{{\rm\bf K}}
\def\C{{\rm\bf BC}}

\def\H{\mathcal{H}}
\def\newH{\mathcal{H}}

\def\R{\mathcal{R}}
\def\demo{\noindent{\bf Proof}\hskip10pt}

\def\qed{\hskip10pt $\Box$\vspace{3mm}}

\tikzset{every picture/.style={line width=0.75pt}} 

\begin{document}
\title{A family of $2$-groups and an associated family of semisymmetric, locally $2$-arc-transitive graphs}
\author
{
Daniel R. Hawtin\\
{\small Faculty of Mathematics, The University of Rijeka}\\[-4pt]
{\small Rijeka, 51000, Croatia}\\[-2pt]
{\small{Email}:\ dan.hawtin@gmail.com}\\[+6pt]
Cheryl E. Praeger\\
{\small Department of Mathematics and Statistics, The University of Western Australia}\\[-4pt]
{\small Crawley, WA 6907, Australia}\\[-2pt]
{\small{Email}:\ cheryl.praeger@uwa.edu.au}\\[+6pt]
Jin-Xin Zhou  \\
{\small School of Mathematics and Statistics, Beijing Jiaotong University}\\[-4pt]
{\small Beijing 100044, P.R. China}\\[-2pt]
{\small{Email}:\ jxzhou@bjtu.edu.cn}
}

\date{\emph{Dedicated to the memory of Zvonimir Janko.}}
\maketitle
\begin{abstract}
A \emph{mixed dihedral group} is a group $H$ with two disjoint subgroups $X$ and $Y$, each elementary abelian of order $2^n$, such that $H$ is generated by $X\cup Y$, and $H/H'\cong X\times Y$. In this paper, for each $n\geq 2$, we construct a mixed dihedral $2$-group $H$ of nilpotency class $3$ and order $2^a$ where $a=(n^3+n^2+4n)/2$, and a corresponding graph $\Sigma$, which is the clique graph of a Cayley graph of $H$. We prove that $\Sigma$ is semisymmetric, that is, $\Aut(\Sigma)$ acts transitively on the edges, but intransitively on the vertices, of $\Sigma$. These graphs are the first known semisymmetric graphs constructed from groups that are not $2$-generated (indeed $H$ requires $2n$ generators). Additionally, we prove that $\Sigma$ is locally $2$-arc-transitive, and is a normal cover of the `basic' locally $2$-arc-transitive graph $\K_{2^n,2^n}$. As such, the construction of this family of graphs contributes to the investigation of normal covers of prime-power order of basic locally $2$-arc-transitive graphs -- the `local' analogue of a question posed by C.~H.~Li.

\bigskip
\noindent{\bf Key words:} semisymmetric, $2$-arc-transitive, edge-transitive, normal cover, Cayley graph \\
\noindent{\bf 2000 Mathematics subject classification:} 05C38, 20B25
\end{abstract}

\section{Introduction}

Many graphs with a lot of symmetry arise from constructions based on groups. These include Cayley graphs, and more generally arc-transitive coset graphs, all of which are vertex-transitive. More recently the Cayley graph construction was extended to a theory of bi-Cayley graphs in \cite{ZF} which led to the construction of the first (to our knowledge) infinite family of semisymmetric graphs based on finite $2$-groups. 

\medskip\noindent
\emph{Semisymmetric graphs.} \quad These are regular graphs (that is, each vertex has the same valency) which are edge-transitive but not vertex-transitive. They have been studied for more than 50 years. In 1967, in what is perhaps the first paper published on the subject, Folkman   gave a method for constructing examples of semisymmetric graphs from abelian groups \cite[Theorem 4]{Fo} and posed a number of questions, notably he asked for \emph{all values of $v, k$ such there exists a semisymmetric graph on $v$ vertices (that is, of order $v$) with valency $k$}. By 2006, all cubic (valency $3$) semisymmetric graphs on up to $768$ vertices had been enumerated, (by Ivanov~\cite{AVIv} for orders up to $28$ in 1987, and the rest by Conder et. al.~\cite{C-M-M-P}). Ming Yao Xu alerted the third author that the list contained no examples with $2$-power order, and it turned out that the theory of bi-Cayley graphs developed in \cite{ZF} could be applied to construct, for each $n\geq 2$, a cubic semisymmetric graph of order $2^{2n+7}$ which is a bi-Cayley graph for a $2$-group $H$ of order $2^{2n+6}$. The group $H$ was $2$-generated with derived quotient $H/H'\cong C_{2^n}\times C_{2^n}$ and $|H'|=2^6$. An additional family of semisymmetric bi-Cayley graphs was given by Conder et. al.~\cite[Example 5.2 and Proposition 5.4]{C-Z-F-Z} in 2020. This time the  graphs had order $4n$ and valency $2k$, with $k$ odd, and were constructed as bi-Cayley graphs for a dihedral group $D_{2n}$ of order $2n$, with the requirement that some element of $\mathbb{Z}_n^*$ has multiplicative order $2k$. Thus although the valencies in this new family were unbounded, the groups used in the construction were still $2$-generated.

Our aim in this paper is to present a new infinite family of semisymmetric graphs based on very different kinds of $2$-groups. For each $n\geq 2$ we construct (see Definitions~\ref{mix-dih},~\ref{def:Hn} and Theorem~\ref{t:semi-sym}) a semisymmetric graph of order $2^{n^2(n+1)/2 + n + 1}$ and valency $2^n$, based on a $2$-group $H$ with $H/H'=C_2^{2n}$ (which implies that $H$ requires $2n$ generators). 
The idea for our construction came from our recent paper \cite{HPZ} where we studied a natural Cayley graph $\Gamma(H)$ for a group $H$ (not necessarily a $2$-group) with disjoint subgroups $X, Y$ such that $X\cong Y\cong C_2^n$, $H=\langle X, Y\rangle$, and $H/H'\cong C_2^{2n}$ (Definition~\ref{mix-dih}). It turns out that, in the case where $\Gamma(H)$ is edge-transitive, its clique graph has many desirable properties \cite[Theorem 1.6]{HPZ}. Here we apply this theory to an explicit family of such groups $H$ with nilpotency class $3$ and show that the associated clique graphs of the Cayley graphs $\Gamma(H)$  are semisymmetric (Theorem~\ref{t:semi-sym}).

\medskip\noindent
\emph{Locally $2$-arc-transitive graphs.}\quad In making this construction we had an additional objective in mind. Our construction produces semisymmetric graphs which are locally $2$-arc-transitive (Theorem~\ref{t:semi-sym}). The $2$-arcs in a graph $\Sigma$ are vertex-triples $(u,v,w)$ such that $\{u,v\}$ and $\{v,w\}$ are edges and $u\ne w$, and $\Sigma$ is said to be {\em locally $(G,2)$-arc-transitive\/} if $G\leq\Aut(\G)$ and for each vertex $u$, the stabiliser $G_u$ is transitive on the $2$-arcs  $(u,v,w)$ starting at $u$.
For a finite connected locally $(G,2)$-arc-transitive graph $\Sigma$, the group $G$ is edge-transitive, and either $G$ is vertex-transitive or $\Sigma$ is bipartite and the two parts of the vertex-bipartition are the $G$ vertex-orbits. Such graphs had been extensively studied since the seminal work of Tutte~\cite{T47} in the 1940s, and more recently the second author (in 1993 \cite[Theorem 4.1]{ONS} for $G$-vertex-transitive graphs, and in 2004 with Giudici and Li~\cite[Theorems 1.2 and 1.3]{GLP-tamc} for $G$-vertex-intransitive graphs) identified a sub-family of `basic' locally $(G,2)$-arc-transitive graphs such that each finite connected locally $(G,2)$-arc-transitive graph is a {\em normal cover} of a basic example (see Section~\ref{s:nquots} for a discussion of these concepts).  This new approach allows effective use of modern permutation group theory and the finite simple group classification to study these graphs.  We note that, if a locally $(G,2)$-arc-transitive, $G$-vertex-intransitive graph is a normal cover of a basic graph $\Sigma_0$, then $\Sigma_0$ may have additional symmetry not inherited from $G$; in particular it may be vertex-transitive. 

The family of `basic' graphs that are covered by the graphs in our construction are the complete bipartite graphs $\K_{2^n,2^n}$, which of course are vertex-transitive. They form one of a small number of families of basic locally $(G,2)$-arc-transitive, $G$-vertex-transitive graphs of prime power order classified in~\cite{L-BLMS-2001}, sharpening the classification arising from \cite{IP, ONS, Bip}  for prime power orders (see Subsection~\ref{s:nquots}). They also arise, as we prove in Theorem~\ref{t:semi-sym}, as basic graphs covered by $G$-vertex-intransitive, locally $(G,2)$-arc-transitive graphs of $2$-power order. It would be good to have an extension of Li's classification (of the vertex-transitive examples) to all basic regular locally $(G,2)$-arc-transitive graphs of prime power order.  We note that there are some examples of vertex-intransitive, basic  locally $(G,2)$-arc-transitive graphs of prime power order that are not regular graphs, for example the stars $K_{1,p^a-1}$ with $G=S_{p^a-1}$, but any regular graphs with these properties will have order a $2$-power. 

\begin{prob}\label{prob-3}
Classify the  basic, regular, locally $(G,2)$-arc-transitive graphs of prime power order. In particular, are there any additional bipartite examples which are not already in Li's classification~
{\rm\cite[Theorem 1.1]{L-BLMS-2001}}? 
\end{prob}
%
Li \cite[pp.130--131]{L-BLMS-2001} was ``inclined to think that non-basic 2-arc-transitive graphs of prime power order would be rare and hard to construct", and posed the problem \cite[Problem]{L-BLMS-2001} of constructing and characterising the normal covers of prime power order of the basic graphs in his classification. We would like to expand this problem to include covers of all basic graphs arising from Problem~\ref{prob-3}. Here, as in \cite{HPZ}, we focus on covers of the graphs $\K_{2^n,2^n}$.


\subsection{The main result}

The general family of groups and graph constructions we will study are specified in Definition~\ref{mix-dih}. One of the graphs is a  {\em Cayley graph} $\Cay(G,S)$ for a group $G$ with respect to an inverse-closed subset $S\subseteq G\setminus\{1\}$ (that is, $s^{-1}\in S$ for all $s\in S$): it is the graph with vertex set $G$ and edge set $\{\{g,sg\}\ :\ g\in G,s\in S\}$.


\begin{defi}\label{mix-dih}
{\em Let $n$ be an integer, $n\geq2$. 

\medskip\noindent
{\rm (a)} If $H$ is a finite group with subgroups $X,Y$ such that $X\cong Y\cong C_2^n$, $H=\lg X, Y\rg$ and $H/H'\cong C_2^{2n}$, where $H'$ is the derived subgroup of $H$, then we say that $H$ is an {\em $n$-dimensional mixed dihedral group relative to $X$ and $Y$}.

    \medskip\noindent
{\rm (b)} For $H, X, Y$ as in part (a),  the graphs $C(H,X,Y)$ and $\Sigma(H,X,Y)$ are defined as follows:
\begin{equation}\label{eq-1}
C(H,X,Y)=\Cay(H,S(X,Y)),\ {\rm with}\ S(X,Y)=(X\cup Y)\setminus\{1\};
\end{equation}
and $\Sigma=\Sigma(H,X,Y)$ is the graph with vertex-set and edge-set given by:
\begin{equation}\label{eq-2}
\begin{array}{l}
V(\Sigma)=\{Xh, Yh: h\in H\},\\
E(\Sigma)=\{\{Xh, Yg\}:  h,g\in H\ \mbox{and}\ Xh\cap Yg\neq\emptyset\}.
\end{array}
\end{equation}
}
\end{defi}


While this construction was used in \cite[Theorem 1.8]{HPZ} to obtain an infinite family of (vertex-transitive) $2$-arc-transitive normal covers of $\K_{2^n, 2^n}$ of order a $2$-power, our interest here is semisymmetric examples.  This is a much more delicate problem, and requires us to analyse a new infinite family of mixed-dihedral $2$-groups $\newH(n)$, which we now define. 

\begin{defi}\label{def:Hn}
{\rm 
Let $n\geq 2$ be an integer and let $X_0=\{x_1,\dots,x_n\}$, $Y_0=\{y_1,\dots,y_n\}$, and consider the group
\[
\newH(n)=\lg X_0\cup Y_0\mid \mathcal{R}\rg
\]
where $\mathcal{R}$ is the following set of relations: for $x,x'\in X_0$, $y,y'\in Y_0$, and  $z,z', z'', z''' \in X_0\cup Y_0$,
\[\begin{array}{l}
z^2=1, [x,x']=[y,y']=1,[x,y]^2=1, [[y,x],y']=1,  [[x,y],z]^2=1, [[[z,z'],z''],z''']=1.
\end{array}\]
}
\end{defi}

It turns out that, for these groups, the graph $\Sigma(\newH(n),X,Y)$ is semisymmetric and is a locally $2$-arc-transitive normal cover of $\K_{2^n, 2^n}$.

\begin{thm}\label{t:semi-sym}
For $n\geq2$, the group $\newH(n)$ in Definition~$\ref{def:Hn}$ is an $n$-dimensional mixed dihedral group of order $2^{n(n^2+n+4)/2}$ relative to $X=\lg X_0\rg$ and $Y=\lg Y_0\rg$, and the graph $\Sigma(\newH(n), X,Y)$ as in $\eqref{eq-2}$ is semisymmetric and locally $2$-arc-transitive, of valency $2^n$ and order $2^{n^2(n+1)/2 +n+1}$.
\end{thm}

\begin{rem}\label{rem:smallestgraph}
{\rm  
 The smallest case in Theorem~\ref{t:semi-sym} is that of $n=2$. Here $|\newH(2)|=2^{10}$, and $\Sigma=\Sigma(\newH(2), X,Y)$ has order $|V(\Sigma)|=2^{9}$ and valency $4$. A computation in GAP \cite{GAP4} shows that 
 $|\Aut(\Sigma)|=2^{15}\cdot 3^5$, considerably larger than the subgroup $A:=\newH(2)\cdot A(\newH(2),X,Y)=\newH(2)\cdot(\GL_2(2)\times\GL_2(2))$ of order $2^{12}\cdot 3^2$ used in the proof of Theorem~\ref{t:semi-sym} (see also Lemma~\ref{lem:prop-mixed-dih} and Theorem~\ref{prop-2groupNew}~(5)). To give some insight into the structure of $\Sigma$ we computed, again using GAP, the  distance diagrams for $\Sigma$ from the vertices $X$ and $Y$. These are shown in Figure~\ref{distdiagSigma}, and demonstrate that $\Sigma$ is locally $3$-distance-transitive. The two diagrams exhibit strikingly different structure from distance $4$ onwards. In fact,  the stabilisers $\Aut(\Sigma)_X$ and $\Aut(\Sigma)_Y$ of the vertices $X$ and $Y$ are non-isomorphic subgroups of order $2^7\cdot 3^5$. Further computations, performed in GAP and Magma \cite{BCP}, show that 
 \begin{enumerate}
     \item $O_2(\Aut(\Sigma)) = \newH(2)'.Y \cong C_2^8$, lies in $\newH(2)$, is faithful and regular on the $\Aut(\Sigma)$-orbit containing the vertex $X$, and has four orbits of length $64$ on the $\Aut(\Sigma)$-orbit containing $Y$. The normal quotient of $\Sigma$ modulo $O_2(\Aut(\Sigma))$ is the `star' $\K_{1,4}$. 
     \item Also $A$ is self-normalising in $\Aut(\Sigma)$, and there are exactly $896$ subgroups of 
 $\Aut(\Sigma)$ of order $256$ that act semiregularly on $V(\Sigma)$ with orbits the two parts of the bipartition. Hence, at least in the case $n=2$, the graph $\Sigma$ is a bi-Cayley graph of some $2$-group.
 \end{enumerate}
 }
\end{rem}

\tikzset{every picture/.style={line width=0.75pt}} 

\begin{figure}
 \[
 \begin{tikzpicture}[x=0.75pt,y=0.75pt,yscale=-1,xscale=1]

\draw   (25,110) .. controls (25,101.72) and (31.72,95) .. (40,95) .. controls (48.28,95) and (55,101.72) .. (55,110) .. controls (55,118.28) and (48.28,125) .. (40,125) .. controls (31.72,125) and (25,118.28) .. (25,110) -- cycle ;
\draw   (375,110) .. controls (375,101.72) and (381.72,95) .. (390,95) .. controls (398.28,95) and (405,101.72) .. (405,110) .. controls (405,118.28) and (398.28,125) .. (390,125) .. controls (381.72,125) and (375,118.28) .. (375,110) -- cycle ;
\draw   (305,110) .. controls (305,101.72) and (311.72,95) .. (320,95) .. controls (328.28,95) and (335,101.72) .. (335,110) .. controls (335,118.28) and (328.28,125) .. (320,125) .. controls (311.72,125) and (305,118.28) .. (305,110) -- cycle ;
\draw   (165,110) .. controls (165,101.72) and (171.72,95) .. (180,95) .. controls (188.28,95) and (195,101.72) .. (195,110) .. controls (195,118.28) and (188.28,125) .. (180,125) .. controls (171.72,125) and (165,118.28) .. (165,110) -- cycle ;
\draw   (95,110) .. controls (95,101.72) and (101.72,95) .. (110,95) .. controls (118.28,95) and (125,101.72) .. (125,110) .. controls (125,118.28) and (118.28,125) .. (110,125) .. controls (101.72,125) and (95,118.28) .. (95,110) -- cycle ;
\draw   (445,110) .. controls (445,101.72) and (451.72,95) .. (460,95) .. controls (468.28,95) and (475,101.72) .. (475,110) .. controls (475,118.28) and (468.28,125) .. (460,125) .. controls (451.72,125) and (445,118.28) .. (445,110) -- cycle ;
\draw    (55,110) -- (95,110) ;
\draw    (125,110) -- (165,110) ;
\draw   (95,205) .. controls (95,196.72) and (101.72,190) .. (110,190) .. controls (118.28,190) and (125,196.72) .. (125,205) .. controls (125,213.28) and (118.28,220) .. (110,220) .. controls (101.72,220) and (95,213.28) .. (95,205) -- cycle ;
\draw   (305,235) .. controls (305,226.72) and (311.72,220) .. (320,220) .. controls (328.28,220) and (335,226.72) .. (335,235) .. controls (335,243.28) and (328.28,250) .. (320,250) .. controls (311.72,250) and (305,243.28) .. (305,235) -- cycle ;
\draw   (305,175) .. controls (305,166.72) and (311.72,160) .. (320,160) .. controls (328.28,160) and (335,166.72) .. (335,175) .. controls (335,183.28) and (328.28,190) .. (320,190) .. controls (311.72,190) and (305,183.28) .. (305,175) -- cycle ;
\draw   (235,205) .. controls (235,196.72) and (241.72,190) .. (250,190) .. controls (258.28,190) and (265,196.72) .. (265,205) .. controls (265,213.28) and (258.28,220) .. (250,220) .. controls (241.72,220) and (235,213.28) .. (235,205) -- cycle ;
\draw   (165,205) .. controls (165,196.72) and (171.72,190) .. (180,190) .. controls (188.28,190) and (195,196.72) .. (195,205) .. controls (195,213.28) and (188.28,220) .. (180,220) .. controls (171.72,220) and (165,213.28) .. (165,205) -- cycle ;
\draw    (125,205) -- (165,205) ;
\draw    (195,205) -- (235,205) ;
\draw    (265,200) -- (305,180) ;
\draw    (265,210) -- (305,230) ;
\draw    (335,230) -- (375,210) ;
\draw   (585,110) .. controls (585,101.72) and (591.72,95) .. (600,95) .. controls (608.28,95) and (615,101.72) .. (615,110) .. controls (615,118.28) and (608.28,125) .. (600,125) .. controls (591.72,125) and (585,118.28) .. (585,110) -- cycle ;
\draw   (515,110) .. controls (515,101.72) and (521.72,95) .. (530,95) .. controls (538.28,95) and (545,101.72) .. (545,110) .. controls (545,118.28) and (538.28,125) .. (530,125) .. controls (521.72,125) and (515,118.28) .. (515,110) -- cycle ;
\draw    (475,110) -- (515,110) ;
\draw    (545,110) -- (585,110) ;
\draw   (25,205) .. controls (25,196.72) and (31.72,190) .. (40,190) .. controls (48.28,190) and (55,196.72) .. (55,205) .. controls (55,213.28) and (48.28,220) .. (40,220) .. controls (31.72,220) and (25,213.28) .. (25,205) -- cycle ;
\draw    (55,205) -- (95,205) ;
\draw    (335,110) -- (375,110) ;
\draw    (405,110) -- (445,110) ;
\draw    (195,110) -- (235,110) ;
\draw    (265,110) -- (305,110) ;
\draw   (235,110) .. controls (235,101.72) and (241.72,95) .. (250,95) .. controls (258.28,95) and (265,101.72) .. (265,110) .. controls (265,118.28) and (258.28,125) .. (250,125) .. controls (241.72,125) and (235,118.28) .. (235,110) -- cycle ;
\draw   (445,235) .. controls (445,226.72) and (451.72,220) .. (460,220) .. controls (468.28,220) and (475,226.72) .. (475,235) .. controls (475,243.28) and (468.28,250) .. (460,250) .. controls (451.72,250) and (445,243.28) .. (445,235) -- cycle ;
\draw   (445,175) .. controls (445,166.72) and (451.72,160) .. (460,160) .. controls (468.28,160) and (475,166.72) .. (475,175) .. controls (475,183.28) and (468.28,190) .. (460,190) .. controls (451.72,190) and (445,183.28) .. (445,175) -- cycle ;
\draw   (375,205) .. controls (375,196.72) and (381.72,190) .. (390,190) .. controls (398.28,190) and (405,196.72) .. (405,205) .. controls (405,213.28) and (398.28,220) .. (390,220) .. controls (381.72,220) and (375,213.28) .. (375,205) -- cycle ;
\draw    (405,200) -- (445,180) ;
\draw    (405,210) -- (445,230) ;
\draw    (475,230) -- (515,210) ;
\draw   (585,205) .. controls (585,196.72) and (591.72,190) .. (600,190) .. controls (608.28,190) and (615,196.72) .. (615,205) .. controls (615,213.28) and (608.28,220) .. (600,220) .. controls (591.72,220) and (585,213.28) .. (585,205) -- cycle ;
\draw   (515,205) .. controls (515,196.72) and (521.72,190) .. (530,190) .. controls (538.28,190) and (545,196.72) .. (545,205) .. controls (545,213.28) and (538.28,220) .. (530,220) .. controls (521.72,220) and (515,213.28) .. (515,205) -- cycle ;
\draw    (545,205) -- (585,205) ;

\draw (40,110) node  [font=\small]  {$1$};
\draw (110,110) node  [font=\small]  {$4$};
\draw (57,91.4) node [anchor=north west][inner sep=0.75pt]  [font=\footnotesize]  {$4$};
\draw (86,91.4) node [anchor=north west][inner sep=0.75pt]  [font=\footnotesize]  {$1$};
\draw (126,91.4) node [anchor=north west][inner sep=0.75pt]  [font=\footnotesize]  {$3$};
\draw (320,235) node  [font=\small]  {$72$};
\draw (320,175) node  [font=\small]  {$9$};
\draw (250,205) node  [font=\small]  {$36$};
\draw (180,205) node  [font=\small]  {$12$};
\draw (156,186.4) node [anchor=north west][inner sep=0.75pt]  [font=\footnotesize]  {$1$};
\draw (196,186.4) node [anchor=north west][inner sep=0.75pt]  [font=\footnotesize]  {$3$};
\draw (226,186.4) node [anchor=north west][inner sep=0.75pt]  [font=\footnotesize]  {$1$};
\draw (266,181.4) node [anchor=north west][inner sep=0.75pt]  [font=\footnotesize]  {$1$};
\draw (266,216.4) node [anchor=north west][inner sep=0.75pt]  [font=\footnotesize]  {$2$};
\draw (373,223) node  [font=\footnotesize]  {$2$};
\draw (343,237) node  [font=\footnotesize]  {$3$};
\draw (298,237) node  [font=\footnotesize]  {$1$};
\draw (298,172) node  [font=\footnotesize]  {$4$};
\draw (600,110) node  [font=\small]  {$81$};
\draw (530,110) node  [font=\small]  {$108$};
\draw (506,91.4) node [anchor=north west][inner sep=0.75pt]  [font=\footnotesize]  {$1$};
\draw (546,91.4) node [anchor=north west][inner sep=0.75pt]  [font=\footnotesize]  {$3$};
\draw (576,91.4) node [anchor=north west][inner sep=0.75pt]  [font=\footnotesize]  {$4$};
\draw (40,205) node  [font=\small]  {$1$};
\draw (57,186.4) node [anchor=north west][inner sep=0.75pt]  [font=\footnotesize]  {$4$};
\draw (180,110) node  [font=\small]  {$12$};
\draw (156,91.4) node [anchor=north west][inner sep=0.75pt]  [font=\footnotesize]  {$1$};
\draw (196,91.4) node [anchor=north west][inner sep=0.75pt]  [font=\footnotesize]  {$3$};
\draw (250,110) node  [font=\small]  {$36$};
\draw (226,91.4) node [anchor=north west][inner sep=0.75pt]  [font=\footnotesize]  {$1$};
\draw (266,91.4) node [anchor=north west][inner sep=0.75pt]  [font=\footnotesize]  {$3$};
\draw (320,110) node  [font=\small]  {$54$};
\draw (296,91.4) node [anchor=north west][inner sep=0.75pt]  [font=\footnotesize]  {$2$};
\draw (336,91.4) node [anchor=north west][inner sep=0.75pt]  [font=\footnotesize]  {$2$};
\draw (390,110) node  [font=\small]  {$108$};
\draw (366,91.4) node [anchor=north west][inner sep=0.75pt]  [font=\footnotesize]  {$1$};
\draw (406,91.4) node [anchor=north west][inner sep=0.75pt]  [font=\footnotesize]  {$3$};
\draw (460,110) node  [font=\small]  {$108$};
\draw (436,91.4) node [anchor=north west][inner sep=0.75pt]  [font=\footnotesize]  {$3$};
\draw (476,91.4) node [anchor=north west][inner sep=0.75pt]  [font=\footnotesize]  {$1$};
\draw (460,235) node  [font=\small]  {$108$};
\draw (460,175) node  [font=\small]  {$27$};
\draw (390,205) node  [font=\small]  {$108$};
\draw (408,188) node  [font=\footnotesize]  {$1$};
\draw (408,222) node  [font=\footnotesize]  {$1$};
\draw (506,216.4) node [anchor=north west][inner sep=0.75pt]  [font=\footnotesize]  {$3$};
\draw (483,237) node  [font=\footnotesize]  {$3$};
\draw (438,237) node  [font=\footnotesize]  {$1$};
\draw (432,165.4) node [anchor=north west][inner sep=0.75pt]  [font=\footnotesize]  {$4$};
\draw (600,205) node  [font=\small]  {$27$};
\draw (530,205) node  [font=\small]  {$108$};
\draw (546,186.4) node [anchor=north west][inner sep=0.75pt]  [font=\footnotesize]  {$1$};
\draw (576,186.4) node [anchor=north west][inner sep=0.75pt]  [font=\footnotesize]  {$4$};
\draw (110,205) node  [font=\small]  {$4$};
\draw (87,186.4) node [anchor=north west][inner sep=0.75pt]  [font=\footnotesize]  {$1$};
\draw (126,186.4) node [anchor=north west][inner sep=0.75pt]  [font=\footnotesize]  {$3$};

\end{tikzpicture}
 \]
 \caption{Distance diagrams at the vertices $X$ (upper) and $Y$ (lower) for the smallest graph $\Sigma=\Sigma(\newH(2),X,Y)$ in Theorem~\ref{t:semi-sym}. It has $512$ vertices and valency $4$. Here each node represents an orbit of the stabiliser of the relevant vertex ($X$ or $Y$) in the full automorphism group of $\Sigma$. Computations performed in GAP \cite{GAP4}.  
 }
 \label{distdiagSigma}
\end{figure}
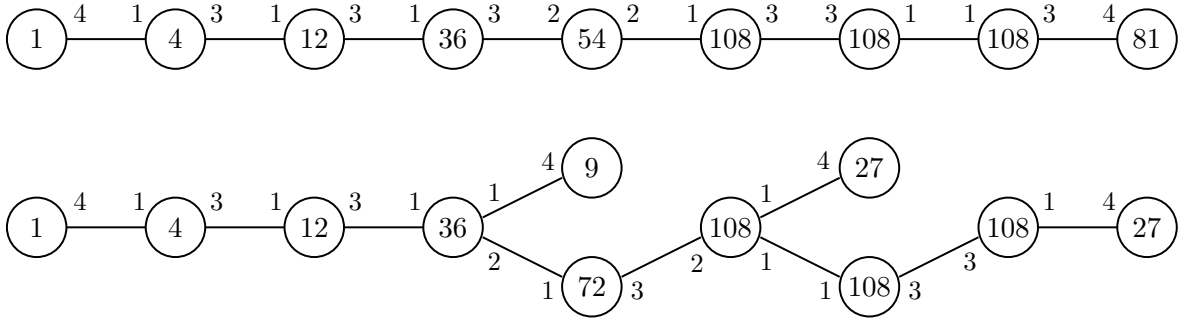



    
This paper is organised as follows. In Section~\ref{prelimSect}, we outline the notation used in the paper and give several preliminary results, including Lemma~\ref{lem:prop-mixed-dih} which summarises several properties of the graphs $C(H,X,Y)$ and $\Sigma(H,X,Y)$ which we will need, and which were proved in \cite{HPZ}. In Section~3, we investigate the structure of the group $\H(n)$, and in Section~4 we prove Theorem~\ref{t:semi-sym}.

\section{Notation and preliminary results for graphs}\label{prelimSect}


All graphs we consider are finite, connected, simple and undirected. 
Let $\G$ be a graph. Denote by $V(\Gamma)$, $E(\Gamma)$ and $\Aut(\Gamma)$ the vertex set, edge set, and full automorphism group of $\G$, respectively. For $v\in V(\G)$, let $\G(v)$ denote the set of vertices adjacent to $v$. A graph $\G$ is said to be {\em regular} if there exists an integer $k$ such that $|\G(v)|=k$ for all vertices $v\in V(\G)$. A graph $\G$ is bipartite if $E(\G)\ne\emptyset$ and $V(\G)$ is of the form $\Delta\cup\Delta'$ such that each edge consists of one vertex from $\Delta$ and one vertex from $\Delta'$. If $\G$ is connected then this vertex partition is uniquely determined and the two parts $\Delta, \Delta'$ are often called the \emph{biparts} of $\G$.

For a graph $\G$, let $G\leq\Aut(\G)$. For $v\in V(\G)$, let $G_v=\{g\in G\ :\ v^g=v\}$, the stabiliser of $v$ in $G$. We say that $\G$ is {\em $G$-vertex-transitive\/} or {\em $G$-edge-transitive\/} if $G$ is transitive on $V(\G)$ or $E(\G)$, respectively, and that $\G$ is {\em $G$-semisymmetric\/} if $\G$ is regular and $G$-edge-transitive but not $G$-vertex-transitive. When $G=\Aut(\G)$, a $G$-vertex-transitive, $G$-edge-transitive or $G$-semisymmetric graph $\G$ is simply called {\em vertex-transitive\/}, {\em edge-transitive\/} or {\em semisymmetric\/}, respectively. The $2$-arcs in a graph $\G$ are vertex-triples $(u,v,w)$ such that $\{u,v\},\{v,w\}\in E(\G)$ and $u\ne w$.  A  graph $\G$ is said to be {\em locally $(G,2)$-arc-transitive\/} if $G\leq\Aut(\G)$ and, for each $u\in V(\G)$, $G_u$  is transitive on the $2$-arcs  $(u,v,w)$ starting at $u$, or equivalently, see \cite[Lemma 3.2]{GLP-tamc}, $G_u$ is $2$-transitive on the set $\G(u)$.
Similarly, when $G=\Aut(\G)$, a locally $(G,2)$-arc-transitive graph $\G$ is simply called {\em locally $2$-arc-transitive\/}.  There is a considerable body of literature on locally 2-arc-transitive graphs, see for example \cite{FanLLP2013,GLP-tamc,L-TAMS-2006,Li-Pan,ONS}. 

\subsection{Normal quotients and normal covers of graphs}\label{s:nquots}
 The normal quotient method for investigating vertex- or edge-transitive graphs proceeds as follows. 
 Assume that $G\leq\Aut(\G)$ is such that $\G$ is $G$-vertex-transitive or $G$-edge-transitive. Let $N$ be a normal subgroup of $G$ such that $N$ is intransitive on $V(\G)$. The {\em $N$-normal quotient graph\/} of $\G$ is defined as the graph $\G_N$ with vertices the $N$-orbits in $V(\G)$ and with two distinct $N$-orbits adjacent if there exists an edge in $\G$ consisting of one vertex from each of these orbits. If $\G$ is regular, and if $\G_N$ and $\G$ have the same valency, then we say that $\G$ is an {\em $N$-normal cover} of $\G_N$.
 

If $\G$ is connected and is a regular, locally $(G,2)$-arc-transitive graph, and if $N$ is intransitive on each $G$-vertex-orbit, then  (by \cite[Theorem 4.1]{ONS} and \cite[Lemma 5.1]{GLP-tamc}) also $\G_N$ is a connected, regular locally $(G/N,2)$-arc-transitive graph, $\G$ is an $N$-normal cover of $\G_N$, and $N$ is semiregular on $V(\G)$, that is, each $N$-orbit has size $|N|$. Such a graph $\G$ is said to be {\em basic} (or sometimes $G$-basic, to emphasise dependence on $G$) if there is no suitable normal subgroup $N$ to make such a reduction; that is, if each nontrivial normal subgroup $N$ of $G$ is transitive on at least one $G$-orbit, forcing the quotient $\G_N$ to be degenerate, namely either $\K_1$ or a star $\K_{1,k}$ for some $k\geq1$. In some cases the $G$-basic graphs can be determined:  
Li's classification in \cite{L-BLMS-2001} shows that the $G$-basic graphs of prime power order, in the case where $G$ is vertex-transitive are one of: $\K_{2^n,2^n}$ (the complete bipartite graph), $\K_{p^m}$ (the complete graph), $\K_{2^n,2^n}-2^n\K_2$ (the graph obtained by deleting a $1$-factor from $\K_{2^n,2^n}$) or a primitive or biprimitive `affine graph' (by which, see \cite[p.130]{L-BLMS-2001}, Li meant the graphs  in the classification by Ivanov and the second author in \cite[Table 1]{IP}).

\subsection{Cliques, clique graphs and line graphs}

A \emph{clique} of a graph $\G$ is a subset  $U\subseteq
V(\G)$ such that every pair of vertices in $U$ forms an
edge of $\G$. A clique $U$ is \emph{maximal} if no subset
of $V(\G)$  properly containing $U$ is a clique. The \emph{clique graph}  of $\G$ is defined as the graph $\Sigma(\G)$ with vertices the maximal cliques of $\G$ such that two distinct maximal cliques are adjacent in $\Sigma(\G)$ if and only if their intersection is non-empty. Similarly the \emph{line graph}  of $\G$ is defined as the graph $\mathcal{L}(\G)$ with vertex set $E(\G)$ such that two distinct edges $e,e'\in E(\G)$ are adjacent in $\mathcal{L}(\G)$ if and only if $e\cap e'\ne\emptyset$.

\subsection{Cayley graphs and bi-Cayley graphs}\label{sub:cay}

A group $G$ of permutations of a set $V(\G)$ is called \emph{regular} if it is transitive, and some (and hence all) stabilisers $G_v$ are trivial. (It is unfortunate that this conflicts with the usage of `regular' as defined above for graphs.) More generally $G$ is called \emph{semiregular} if the stabiliser $G_v=1$ for all $v\in V(\G)$. So $G$ is regular if and only if it is semiregular and transitive.

Let $\G=\Cay(G,S)$ be a Cayley graph on $G$ with respect to $S$. For any $g\in G$ define
\[
R(g): x\mapsto xg\ \mbox{for $x\in G$ and set $R(G)=\{R(g)\ :\ g\in G\}$.}
\]
Then $R(G)$ is a regular permutation group on $V(\G)$ (see, for example \cite[Lemma 3.7]{PS}) and is a subgroup of $\Aut(\Cay(G,S))$ (as $R(g)$ maps each edge $\{x,sx\}$ to an edge $\{xg,sxg\}$). For briefness, we shall identify $R(G)$ with $G$ in the following.
Let
\[
\Aut(G,S)=\{\a\in\Aut(G): S^\a=S\}.
\]
It was proved by Godsil~\cite{Godsil-1981} that the normaliser of $G$ in $\Aut(\Cay(G,S))$ is $G: \Aut(G,S)$. 

Cayley graphs are precisely those graphs $\Gamma$ for which $\Aut(\G)$ has a subgroup that is regular on $V(\G)$. For this reason we say that a  graph $\Gamma$ is a \emph{bi-Cayley graph} if $\Aut(\G)$ has a subgroup $H$ which is semiregular on $V(\G)$ with two orbits. The following result from \cite{HPZ} will be useful.

\begin{lem}{\rm\cite[Lemma~2.6]{HPZ}}\label{sufficient}
Let $\G$ be a connected $(G,2)$-arc-transitive graph, and let $u\in V(\G)$. Suppose that $\G$ is an $N$-normal cover of $\K_{2^n,2^n}$, for some normal $2$-subgroup $N$ of $G$. Then $\G$ is bipartite,  and one of the following holds:
\begin{enumerate}
  \item [{\rm (1)}] $\G$ is a normal Cayley graph of a $2$-group;
  \item [{\rm (2)}] $\G$ is a bi-Cayley graph of a $2$-group $H$ such that $G\leq N_{\Aut(\Gamma)}(H)$;
  \item [{\rm (3)}] $N\unlhd\Aut(\G)$.
\end{enumerate}
Moreover if the stabiliser $G_u$ acts unfaithfully on $\G(u)$, then part $(3)$ holds.
\end{lem}

The next result is developed from \cite[Theorem 1.1]{LMP2009} for the case of locally primitive bipartite graphs $\G$. For $G\leq \Aut(\G)$, we denote by  $G^+$, the subgroup of $G$ (of index at most $2$) that fixes both biparts of $V(\G)$ setwise. The proof uses the following concept: a permutation group $G\leq {\rm Sym}(V)$ is \emph{biquasiprimitive} if each nontrivial normal subgroup has at most two orbits in $V$, and there exists such a subgroup having two orbits. 

\begin{lem}\label{reducetoK2n2n}
Let $\G$ be a connected bipartite graph of order $2^m$ and valency $2^n$ with $m>n$. Assume that $\G$ is vertex-transitive and locally primitive, and that $\Aut(\G)_u$ is not faithful on $\G(u)$, for some $u\in V(\G)$. Then $n\geq2$, and there exists $N\unlhd\Aut(\G)$ such that $N$ is a $2$-group, $N\leq \Aut(\G)^+$ and  is semiregular on $V(\G)$,  $\G$ is an $N$-normal cover of $\G_N$, and $\G_N\cong\K_{2^n,2^n}$. 
\end{lem}


\f\demo Let $X=\Aut(\G)$, and note that $n\geq2$, since if $n=1$ then $\G$ is a cycle and $X_u\cong C_2$ is faithful on $\G(u)$. Let $N\unlhd X$ be maximal subject to the condition that $N$ has at
least three orbits on $V(\G)$ (possibly $N=1$), let $\overline{X}=X/N$, and let $X^+$ be the normal subgroup of index $2$ in $X$ which fixes both parts of the bipartition of $\G$. By \cite[Lemma 1.6]{Imp}, or see \cite[Theorem 1.3]{LPVZ}, $N$ is semiregular on $V(\G)$, and hence $N$ is a $2$-group since $|V(\G)|=2^m$. Also, by \cite[Lemma 3.1]{LPVZ}, the quotient $\G_N$ is bipartite and $N\leq X^+$, so $X^+/N$ is an index $2$ normal subgroup of $\overline{X}$ with orbits in $V(G_N)$ the two biparts of $\G_N$. Moreover, by the maximality of $N$, each nontrivial normal subgroup of $\overline{X}$ has at most two orbits in $V(\G_N)$,   and hence $\overline{X}$ is bi-quasiprimitive on $V(\G_N)$.  By \cite[Theorem 1.3]{LPVZ}, $\G_N$ is $\overline{X}$-locally primitive and $\G$ is an $N$-normal cover of $\G_N$. Thus $\G_N$ also has valency $2^n$, and $\G_N$ has order $2^k = |V(\G)|/|N|$ for some $k>n\geq2$. If $\G_N\cong\K_{2^n,2^n}$ then the result holds, so we may assume that  $\G_N\ncong\K_{2^n,2^n}$.



  Then by the third paragraph in the proof of \cite[Theorem~1.1]{LMP2009} (and using \cite[Lemmas~3.3 and 4.2]{LMP2009}), $X$ has a subgroup $G$ satisfying $N<G\leq X^+$ such that $G$ is faithful and regular on each of the biparts of $V(\G)$, and moreover $\G$ is a bi-Cayley graph of $G$ of the form $\G= \Upsilon \times\K_2$ (the direct product of $\Upsilon$ and $\K_2$), where $\Upsilon=\Cay(G,S)$ is a Cayley graph of $G$ with respect to some subset $S$ of $G$.
  Further, $\G_N$ and $\overline{X}$ satisfy \cite[Lemma~3.3(ii) or (iii)]{LMP2009}. If \cite[Lemma~3.3(ii)]{LMP2009} holds then by the second last paragraph  in the proof of \cite[Theorem~1.1]{LMP2009}, $\G$ is a normal Cayley graph of $G\times C_2$. However in this case $X_u$ is faithful on $\G(u)$, which is a contradiction. Thus \cite[Lemma~3.3(iii)]{LMP2009} holds, and then the last paragraph  in the proof of \cite[Theorem~1.1]{LMP2009} shows (since $|V(\G)|=2^m$) that the quotient $\Upsilon_N=\K_{2^r}^{\times \ell}$ (a direct product of $\ell$ copies of $\K_{2^r}$) of valency $(2^r-1)^\ell$.   However in this case, the three graphs $\G$, $\Upsilon$ and $\Upsilon_N$ have the same valency $2^n\geq4$,  which is a contradiction.\hfill\qed
  

\subsection{Mixed dihedrants and their clique graphs}

We record the properties we will need of the graphs from Definition~\ref{mix-dih}.

\begin{lem}{\rm\cite[Lemmas~4.1-4.2]{HPZ}}\label{lem:prop-mixed-dih}
Let $H$ be an $n$-dimensional mixed dihedral group relative to $X$ and $Y$ with $|X|=|Y|=2^n\geq4$, and let $C(H,X,Y)$ and $\Sigma(H,X,Y)$ be the graphs  defined in Definition~\ref{mix-dih}. let $\Sigma =\Sigma(H,X,Y)$, and $G=H:A(H,X,Y)$, where $A(H,X,Y)$ is the setwise stabiliser in $\Aut(H)$ of $X\cup Y$. Then the following hold.
\begin{enumerate}
\item [{\rm (1)}] $\Sigma(H,X,Y)$ is the clique graph of $C(H,X,Y)$.

\item [{\rm (2)}] The map $\varphi:z\to \{Xz,Yz\}$, for $z\in H$, is a bijection $\varphi:H\to E(\Sigma(H,X,Y))$, and induces a graph isomorphism from  $C(H,X,Y)$ to the line graph $\mathcal{L}(\Sigma(H,X,Y))$ of $\Sigma(H,X, Y)$.

\item [{\rm (3)}] $\Aut(C(H,X,Y))=\Aut(\Sigma(H,X,Y))= \Aut(\mathcal{L}(\Sigma(H,X,Y)))$.

\item [{\rm (4)}] The group $G$ acts as a subgroup of automorphisms on $\Sigma$ as follows, for $h,z\in H, \sigma\in A(H,X,Y)$, and $\varphi:H\to E(\Sigma)$ as in part {\rm (2)}:
\begin{align*}
    \mbox{Vertex action:} && h:Xz\to Xzh,\ &  Yz\to Yzh && \sigma:Xz\to X^\sigma z^\sigma,\ \  Yz\to Y^\sigma z^\sigma\\
    \mbox{Edge action:} &&  h:\varphi(z)\to \varphi(zh) &&& \sigma:\varphi(z)\to \varphi(z^\sigma)
\end{align*}
The subgroup $H$ acts regularly on $E(\Sigma)$ and has two orbits on $V(\Sigma)$. In particular, this $G$-action is edge-transitive.

\item [{\rm (5)}] The $H'$-normal quotient graph $\Sigma_{H'}$ of $\Sigma$ is isomorphic to $\K_{2^n,2^n}$ and admits $G/H'$ as an edge-transitive group of automorphisms. Moreover, $\Sigma$ is an $H'$-normal cover of $\K_{2^n,2^n}$.

\item [{\rm (6)}] $A(H, X, Y))\cong A(H,X,Y)^{X\cup Y}\leq (\Aut(X)\times\Aut(Y)): C_2\cong(\GL(n,2)\times\GL(n,2)): C_2$ where the $C_2$ interchanges $X$ and $Y$.

\end{enumerate}
\end{lem}





\section{Notation and preliminary results for groups}\label{s:prelimgp}

For a positive integer $n$,  $C_n$ denotes a cyclic group of order $n$, and $D_{2n}$ denotes a dihedral group of order $2n$.
For a group $G$, we denote by $1$, $Z(G)$, $\Phi(G)$, $G'$, $\soc(G)$ and $\Aut(G)$, the identity element, the centre, the Frattini subgroup, the derived subgroup, the socle and the automorphism group of $G$, respectively. For a subgroup $H$ of a group $G$,
denote by $C_G(H)$ the centraliser of $H$ in $G$ and by $N_G(H)$ the
normaliser of $H$ in $G$. For elements $a,b$ of $G$, the {\em commutator} of $a,b$ is defined as $[a,b]=a^{-1}b^{-1}ab$.
If $X,Y\subseteq G$, then $[X,Y]$ denotes the subgroup generated by all the commutators $[x,y]$ with $x\in X$ and $y\in Y$. We will need the following result concerning $p$-groups.

\begin{lem}\label{basis}{\rm \cite[5.3.2]{Robinson}}
Let $G$ be a finite $p$-group, for some prime $p$, and let $p^r=|G:\Phi(G)|$. Then, $\Phi(G)=G'G^p$, where $G^p=\lg g^p\ |\ g\in G\rg$.
Moreover, every generating set for $G$ has an $r$-element subset which also generates $G$, and in particular, $G/\Phi(G)\cong C_p^r$.
\end{lem}

\subsection{Some results on commutators in arbitrary groups}
We first cite the so-called Witt-Hall formula.

\begin{lem}{\rm \cite[10.2.1.4]{Hall-book} or \cite[5.1.5(iv)]{Robinson}}\label{Witt-formula}
Let $x,y,z$ be elements of a group $G$. Then
\[
[[x,y^{-1}],z]^y\cdot[[y,z^{-1}],x]^z\cdot [[z,x^{-1}],y]^x=1.\]
\end{lem}

Using the Witt-Hall formula above, we obtain the following lemma.

\begin{lem}\label{Witt-formula-for-meta-abelian}
Let $a,b,c$ be elements of a group $G$ such that $G'$ is abelian. Then
\[ [[a,b],c]\cdot[[b,c],a]\cdot[[c,a],b]=1.\]
\end{lem}

\f\demo Let $a,b,c\in G$. By a direct computation or \cite[5.1.5(iii)]{Robinson}, we have $[b,a^{-1}]^a=[b,a]^{-1}=[a,b]$. Hence  
$[[b,a^{-1}],c]^a=[[b,a^{-1}]^a,c^a]=[[a,b],c[c,a]].$

Again (by direct computation or \cite[5.1.5(ii)]{Robinson}),  $[[a,b],c[c,a]]=[[a,b],[c,a]]\cdot [[a,b],c]^{[c,a]}$. Since $G'$ is abelian, this becomes $[[a,b],c[c,a]]=[[a,b],c]$. Consequently, $[[b,a^{-1}],c]^a=[[a,b],c].$
Similarly, 
$ [[a,c^{-1}],b]^c=[[c,a],b] \quad\text{and}\quad
 [[c,b^{-1}],a]^b=[[b,c],a].$
 
Now applying the Witt--Hall formula from Lemma~\ref{Witt-formula} with $y=a,x=b,z=c$, yields the asserted formula: $[[a,b],c]\cdot[[b,c],a]\cdot[[c,a],b]=1.$
%
\hfill\qed


\subsection{Basic commutators in free groups}
In this section, we recall some theory for commutators of a free group, following Marshall Hall's \cite[Charpter 11]{Hall-book}.

\subsubsection{Formal commutators in free groups}\label{s:formalcomms}
Let $F$ be the free group on the ordered alphabet $A=\{a_1, a_2, \ldots, a_{r}\}$, where $r\geq 1$. For $j\geq 1$, the {\em formal commutator} $c_j$ of $F$, and its \emph{weight} $w(c_j)$, are defined by the rules:
\begin{enumerate}
  \item [{\rm (1)}] For $j=1,2,\ldots,r$, $c_j=a_j$, and these are the commutators of weight $1$; \emph{i.e.}, $w(a_j)=1$.
  \item [{\rm (2)}] If $c_i$ and $c_j$ are (formal) commutators, then $[c_i,c_j]$ is a (formal) commutator, say  $c_k$, and $w(c_k)=w(c_i)+w(c_j)$.
\end{enumerate}
\subsubsection{Basic commutators of weight $\ell$ in free groups}\label{s:basiccomms}
Let $F$ be the free group on the ordered alphabet $A=\{a_1, a_2, \ldots, a_{r}\}$, where $r\geq 1$.
For each positive integer $\ell$, we define as follows the set $\C_\ell$  of {\em basic commutators} of $F$ of weight $\ell$, together with a total ordering on $\cup_{u\geq 1}\C_u$:
\begin{enumerate}
  \item [{\rm (1)}]  $\C_1=\{a_1,\dots,a_r\}$, and we choose the ordering $a_1<a_2<\dots<a_r$.

Let $\ell>1$, and assume inductively that $\cup_{1\leq u<\ell} \C_u$ has been defined and ordered.

  \item [{\rm (2)}] Then the set $\C_\ell$ consists of all the commutators $[c_i, c_j]$ that satisfy the following three conditions:
      \begin{enumerate}
        \item [{\rm (a)}] $c_i, c_j\in \cup_{1\leq u<\ell} \C_u$ with $\ell=w(c_i)+w(c_j)$;
        \item [{\rm (b)}]  $c_j<c_i$;
        \item [{\rm (c)}]  If $c_i=[c_s,c_t]$, where $c_s,c_t\in \cup_{1\leq u<\ell} \C_u$, then $c_t\leq  c_j$. 
      \end{enumerate}
  \item [{\rm (3)}] The ordering on $\cup_{1\leq u<\ell} \C_u$ is extended to $\cup_{1\leq u\leq \ell} \C_u$ as follows: we choose an arbitrary order on the set $\C_\ell$, and if $c_i\in\C_\ell$ and $c_j\in \cup_{1\leq u<\ell} \C_u$, we define $c_j<c_i$. 
\end{enumerate}
We next record the nature and sizes of the two sets $\C_2$ and $\C_3$, using the following arithmetic facts.

\begin{lem}\label{lem:arith}
 Let $n\in\mathbb{N}$ with $n\geq 2$.  Then
 \begin{enumerate}
     \item [{\rm (a)}]  $|\{(i,j): 1\leq j<i\leq n\}| = \frac{n(n-1)}{2}$;
     \item [{\rm (b)}]  $|\{(i,j,k): 1\leq j<i\leq n, \ \mbox{and}\ j< k\leq n\}| = \frac{n(n-1)(2n-1)}{6}
     = \frac{n^3}{3} - \frac{n^2}{2} + \frac{n}{6}$;
      \item [{\rm (c)}]  $|\{(i,j,k): 1\leq j<i\leq n, \ \mbox{and}\ j< k\leq n, k\ne i\}| = \frac{n(n-1)(n-2)}{3} = \frac{n^3}{3} - n^2 + \frac{2n}{3}$; 
     \item [{\rm (d)}]  $|\{(i,j,k): 1\leq j<i\leq n, \ \mbox{and}\ j\leq  k\leq n\}| =  \frac{n^3-n}{3}$.
 \end{enumerate}
\end{lem}

\f\demo 
(a) For each $j$ such that $1\leq j< n$, there are precisely $n-j$ choices for $i$, and hence the number of these pair $(i,j)$ is $\sum_{j=1}^{n-1} (n-j)=\sum_{\ell=1}^{n-1}\ell=\frac{n(n-1)}{2}$.

(b) For a fixed $j$ such that $1\leq j< n$, the number of choices of $i, k$ such that $j<i\leq n$ and $j<k\leq n$ is $(n-j)^2$, yielding a total of  $(n-1)^2+(n-2)^2+\cdots+1 = \frac{(n-1)n(2n-1)}{6}= \frac{n^3}{3} - \frac{n^2}{2} + \frac{n}{6}$ triples $(i,j,k)$ with the required constraints.

(c) The number of these triples is equal to the number of triples in part (b) minus the number of triples $(i,j,i)$ with $1\leq j<i\leq n$, which is $\frac{n(n-1)}{2}$ by part (a).

(d) The number of these triples is equal to the number of triples in part (b) plus the number of triples $(i,j,j)$ with $1\leq j<i\leq n$, which is $\frac{n(n-1)}{2}$ by part (a).
\hfill\qed

\begin{lem}\label{lem-basic-com1}
Let $F$ be the free group on the ordered alphabet $A=\{a_1,a_2,\ldots,a_{r}\}$, where $r\geq1$.
Then 
\begin{align*}
    \C_2&=\{[a_i, a_j]: 1\leq j<i\leq r\}, &\mbox{of size}\ & |\C_2|=\frac{r(r-1)}{2},\\
    \C_3&=\{[[a_i, a_j], a_k] : 1\leq j<i\leq r,  \mbox{and} \ j\leq k\leq r\} &\mbox{of size}\  & |\C_3|=\frac{r^3-r}{3}.
\end{align*}
\end{lem}

\f\demo The set $\C_2$ is as claimed (by conditions (2)(a) and (2)(b) in the definition of $\C_\ell$ with $\ell=2$), and $|\C_2|=\frac{r(r-1)}{2}$ as it is in bijection with the $2$-subsets of $\{1,\ldots,r\}$. Next, each element of $\C_3$ has the form $[c_t,c_u]$ with $w(c_t)=2, w(c_u)=1$ (since $c_u<c_t$), and so $c_u=a_k$ for some $k\in[1,r]$ (by condition 1) and $c_t=[a_i,a_j]$  with $1\leq j<i\leq r$ (as we have just seen) and $j\leq k$ (by condition 2(c)). Thus $\C_3$ is as claimed.
By Lemma~\ref{lem:arith}(c), the number of elements $[[a_i, a_j], a_k]\in \C_3$ with $1\leq j<i\leq r,  \mbox{and} \ j\leq k\leq r$ is $\frac{r^3-r}{3}=|\C_3|$. 
\hfill\qed

The following result is known as the Basis Theorem where, following \cite{Hall-book},  for each $k$ we denote the $k$-th term of the lower central series of a free group $F$ by $F_k$. Note in particular that $F_2=F'$, the derived subgroup of $F$.

\begin{thm}{\rm\cite[Theorem~11.2.4]{Hall-book}}\label{Hall-basis-theorem1}
Let $F$ be the free group on the ordered alphabet $A=\{a_1, a_2, \ldots, a_{r}\}$, where $r\geq 1$,  let $\ell\geq1$, and suppose that $\cup_{1\leq u\leq \ell}\C_u=\{c_1,\dots,c_t\}$ with $c_1<c_2<\dots<c_t$. Then each $f\in F$    has a unique representation
\[f=c_1^{s_1}c_2^{s_2}\ldots c_t^{s_t}\ \mod F_{\ell+1},\]
for some integers $s_i$ and, modulo $F_{\ell+1}$, $\C_\ell$ forms a basis   for the free Abelian group $F_\ell/F_{\ell+1}$.
%
\end{thm}

We complete this section with three results about the quotient $F/F_4$.

\begin{lem}\label{lem:PropertyofF/F4}
Let $F$, the $F_\ell$ and $\C_\ell$ be as in Theorem~$\ref{Hall-basis-theorem1}$. Then
 \begin{enumerate}
     \item [{\rm (1)}] $F'=F_2=\lg \C_2, F_3\rg$;
     \item [{\rm (2)}] $F_3=\lg \C_3,F_4\rg$ and $F_3/F_4\leq Z(F/F_4)$;
     \item [{\rm (3)}] $(F/F_4)'=F'/F_4$ is abelian.
 \end{enumerate}
\end{lem}

\f\demo
Recall that $F_2$ is the derived subgroup $F'$. By Theorem~\ref{Hall-basis-theorem1}, $\C_2F_3=\{cF_3\mid c\in \C_2\}$ is a basis for the free abelian group
$F'/F_3$, and $\C_3F_4=\{cF_4\mid c\in \C_3\}$ is a basis for the free abelian group
$F_3/F_4$. Thus $F_3=\lg \C_3,F_4\rg$,  $F'=\lg \C_2, F_3\rg=\lg \C_2,\C_3,F_4\rg$ and $F'/F_4=\lg \C_2 F_4,\C_3 F_4\rg$. In particular, part (1) is proved. 

Let $c=[[a_i,a_j],a_k]\in \C_3$. Then, for each $a\in F$, the commutator $[c,a]\in F_4$, and hence $cF_4\in Z(F/F_4)$. Since $F_3=\lg \C_3,F_4\rg$ this implies  that $F_3/F_4\leq Z(F/F_4)$, proving part (2).

Since $F_4<F'$, we have $(F/F_4)'=F'/F_4$. We have shown that $F'/F_4=\lg \C_2 F_4, \C_3 F_4\rg$, and by part (2) each element of $\C_3F_4$ commutes with each element of $\C_2F_4$ and $\C_3F_4$. Thus to prove part (3) it remains to prove that each pair of elements of $\C_2F_4$ commute. So let $[a_i,a_j],c\in \C_2$. Then
\begin{align*}
    ([a_i,a_j]F_4)^{cF_4}&=[a_i^c,a_j^c]F_4=[c^{-1}a_i c,c^{-1} a_j c]F_4=[c^{-1}ca_i[a_i,c],c^{-1}ca_j[a_j,c]]F_4\\
    &=[a_i[a_i,c],a_j[a_j,c]]F_4.
\end{align*}
 For all $k\leq r$, $[c,a_k]\in F_3$ by the definition of $F_3$, and hence $[a_k,c]\in F_3$. 
Therefore, by part (2), $[a_i,c]F_4$ and $[a_j,c]F_4$ are contained in $Z(F/F_4)$. It follows that $([a_i,a_j]F_4)^{cF_4}=[a_i,a_j]F_4$. Thus $[a_i,a_j]F_4$ and $cF_4$ commute, and part (3) is proved.\hfill\qed

\begin{lem}\label{lem:Propertyofcommutators}
Let $F, A$, and the $F_\ell$ be as in Theorem~$\ref{Hall-basis-theorem1}$. Then, for any $a_i, a_j, a_k\in A$,
\begin{enumerate}
  \item [{\rm (1)}]\ $(a_i^2F_4)^{a_jF_4}=a_i^2\cdot [[a_i,a_j],a_i]\cdot[a_i,a_j]^2F_4$,
  \item [{\rm (2)}]\ $([a_i,a_j]^2F_4)^{a_kF_4}=[a_i,a_j]^2\cdot[[a_i,a_j],a_k]^2F_4$,
   \item [{\rm (3)}]\ $[a_i,a_j]=[a_j,a_i]^{-1}$,
  \item [{\rm (4)}]\ $[[a_j,a_i],a_k]F_4=[[a_i,a_j]^{-1},a_k]F_4=([[a_i,a_j],a_k]F_4)^{-1}$,
  \item [{\rm (5)}]\ $[[a_i,a_j],a_k]F_4=([[a_j,a_k],a_i]F_4)^{-1}\cdot[[a_i,a_k],a_j]F_4$.
\end{enumerate}
\end{lem}

\f\demo (1), (2) By Lemma~\ref{lem:PropertyofF/F4}~(2), $[[a_i,a_j],a_i]F_4\in Z(F/F_4)$. We use this fact for the last equalities of the following two computations in $F/F_4$, which prove parts (1) and (2). 
\begin{align*}
(a_i^2F_4)^{a_jF_4}&=a_j^{-1}a_i^2a_jF_4= a_j^{-1}a_ia_ja_i[a_i,a_j]F_4= a_j^{-1}a_ja_i[a_i,a_j]a_i[a_i,a_j]F_4\\
&= a_i^2[a_i,a_j][[a_i,a_j],a_i][a_i,a_j]F_4= a_i^2[[a_i,a_j],a_i][a_i,a_j]^2F_4.\\
([a_i,a_j]^2F_4)^{a_kF_4}&=a_k^{-1}[a_i,a_j]^2a_kF_4=a_k^{-1}[a_i,a_j]a_k[a_i,a_j][[a_i,a_j],a_k]F_4\\
&=a_k^{-1}a_k[a_i,a_j][[a_i,a_j],a_k][a_i,a_j][[a_i,a_j],a_k]F_4=[a_i,a_j]^2[[a_i,a_j],a_k]^2F_4.
\end{align*}

(3) This follows directly from: $[a_j,a_i]=(a_j^{-1}a_i^{-1}a_ja_i)^{-1}=a_i^{-1}a_j^{-1}a_ia_j=[a_i,a_j]$.

(4) By part (3), we have $[[a_j, a_i],a_k]=[[a_i,a_j]^{-1},a_k]$. Then, by a direct computation or \cite[5.1.5(iii)]{Robinson}, we have
$[[a_i,a_j]^{-1},a_k]=([[a_i,a_j],a_k]^{[a_j,a_i]})^{-1}.$ Lemma~\ref{lem:PropertyofF/F4}~(3) then implies that \[[[a_i,a_j]^{-1},a_k]F_4=([[a_i,a_j],a_k]F_4)^{-1},\]
proving part (4).

(5)  By Lemma~\ref{lem:PropertyofF/F4}~(3), $F'/F_4$ is abelian and hence, by Lemma~\ref{Witt-formula-for-meta-abelian}, we have
\[
[[a_i,a_j],a_k]\cdot[[a_j,a_k],a_i]\cdot[[a_k,a_i],a_j]F_4=F_4.
\]
Since, by part (4), $([[a_k,a_i],a_j]F_4)^{-1} = [[a_i,a_k],a_j]F_4$, it follows, again using the fact $F'/F_4$ is abelian, that
\begin{equation*}\label{eq-group2-3}
[[a_i,a_j],a_k]F_4=([[a_j,a_k],a_i]F_4)^{-1}\cdot[[a_i,a_k],a_j]F_4.
\end{equation*}
This proves part (5), completing the proof.\hfill\qed

We end this section by considering a certain  subgroup $K$ of $F$ containing $F_4$.


\begin{lem}\label{lem:Anormalsubgroup}
Let $F, A$, and the $F_\ell$ be as in Theorem~$\ref{Hall-basis-theorem1}$, and let 
\[
K=\lg F_4, a_i^2, [a_i,a_j]^2, [[a_i,a_j],a_k]^2, [[a_i,a_j],a_i]:  i,j,k\in [1,r] \rg.
\]
Then the following hold.
\begin{enumerate}
  \item [{\rm (1)}] $K\unlhd F$, and $K$ is also generated by $F_4\cup B_K$, where $B_K$ is the set
  \[
   \{ a_i^2, [a_i,a_j]^2, [[a_i,a_j],a_k]^2, [[a_i,a_j],a_\ell] : 1\leq j<i\leq r,\ j<k\leq r,\ k\ne i,\ \ell\in \{i,j\}\};
  \]
  
  \item [{\rm (2)}] $F'K/K\cong C_2^{r(r-1)(2r-1)/6}$ and $\{cK: c\in \C_2\cup D_K\}$ is a basis for $F'K/K$, where $D_K=\{[[a_i,a_j],a_k]: 1\leq j<i\leq r,\ j<k\leq r,\ k\ne i\}$, and $\C_2$ is as in Lemma~$\ref{lem-basic-com1}$; 
  \item [{\rm (3)}] $F/F'K\cong C_2^{r}$.
\end{enumerate}
\end{lem}

\demo (1)\quad Since $F=\lg A \rg$, to prove that $K\unlhd F$ it is enough to show that $b^a\in K$ (or equivalently that $(bF_4)^{aF_4}\in K/F_4$), for any $a\in A$ and any element $b$ in the given generating set for $K$. Let $a\in A$. By Lemma~\ref{lem:Propertyofcommutators} (1) and (2) we have  $(a_i^2F_4)^{aF_4}, ([a_i,a_j]^2F_4)^{aF_4}\in K/F_4$ for all $i,j$. The remaining generators, of the form 
$[[a_i,a_j],a_k]^2, [[a_i,a_j],a_i]$, all lie in $F_3$, and by Lemma~\ref{lem:PropertyofF/F4}~(2), 
$F_3/F_4\leq Z(F/F_4)$. It follows  that $(bF_4)^{aF_4}=bF_4\in K/F_4$ for each of these generators also. Thus $K\unlhd F$. 

To prove the second assertion of (1) let $K_0:=\lg F_4, B_K\rg$. We show that each of the given generators for $K$ lies in $K_0$. Each of the generators $a_i^2$ lies in $B_K\subset K_0$. By Lemma~\ref{lem:Propertyofcommutators} (3), $[a_i,a_j]=[a_j,a_i]^{-1}$, and hence each of the generators $[a_i,a_j]^2$ lies in $K_0$. For a generator $x=[[a_i,a_j],a_i]$ with $i,j\in [1,r]$,  if $j=i$ then $x=1\in K_0$, if $j<i$ then $x\in B_K\subset K_0$, while if $j>i$ then $x\in [[a_j,a_i],a_i]^{-1}F_4\subset K_0$ (by Lemma~\ref{lem:Propertyofcommutators}(4)). It remains to consider the generators of the form $x=[[a_i,a_j],a_k]^2$. If $i=j$ then $x=1\in K_0$ so we may assume that $i\ne j$. If $k=i$ then we have just shown that $[[a_i,a_j],a_k]\in K_0$ so $x\in K_0$ and we may assume also that $k\ne i$. If $k=j$ then, 
by Lemma~\ref{lem:Propertyofcommutators} (4), $[[a_i,a_j],a_j]\in [[a_j,a_i],a_j]^{-1}F_4$ which we have shown lies in $K_0$ so again $x\in K_0$. Thus we may assume that $i,j,k$ are pairwise distinct. Let $m=\min\{i,j,k\}$. If $m=j$ then $x\in B_K\subset K_0$; if $m=i$ then by Lemma~\ref{lem:Propertyofcommutators} (4), 
$x\in [[a_j,a_i],a_k]^{-2}F_4 \subset K_0$; and if $m=k$ then by Lemma~\ref{lem:Propertyofcommutators} (4),  modulo $F_4$, $x= [[a_j,a_k],a_i]^{-2}\cdot [[a_i,a_k],a_j]^2$, and each factor lies on $K_0$ so $x\in K_0$. Thus all generators lie in $K_0$ and hence $K_0=K$. 
This completes the proof of part (1).

(2)\quad 
Now $F'K/K\cong F'/(F'\cap K)\cong (F'/F_4)/((F'\cap K)/F_4)$ (note that $F_4\leq K\cap F'$). By Lemma~\ref{lem:PropertyofF/F4}~(3), $F'/F_4$ is abelian, and hence $F'/(F'\cap K)$ is abelian. Moreover, by \cite[Hilfsatz 1.11]{Huppert}, $F'$ is generated by $[a_i,a_j]^g$ for $i,j\in [1,r]$ and $g\in F$. Then  since each $[a_i,a_j]^2\in K$, it follows that $F'K/K\cong C_2^m$, for some $m$, and we need to find $m$. 

Note that, by Lemma~\ref{lem-basic-com1} and the definition of $D_K$, $\C_3$ is the disjoint union $\C_3=D_K\cup S$ where $S=\{[[a_i,a_j], a_\ell]: 1\leq j<i\leq r,\ \ell\in\{ i, j\}\}$ (and note that $S\subset B_K\subset K$). By Lemma~\ref{lem:PropertyofF/F4}, $F'=\lg \C_2,\C_3, F_4\rg=\lg \C_2,D_K, S, F_4\rg$, and since $K$ contains $F_4\cup S$, it follows that $\{xK : x\in \C_2\cup D_K\}$ is a generating set for  $F'K/K\cong C_2^m$. \emph{We will show that it is in fact a basis using Theorem~\ref{Hall-basis-theorem1}.} It is helpful to consider the following subgroup $H$ such that $F_4<H<K$, 
\[
 H=\lg F_4, [a_i,a_j]^2, [[a_i,a_j],a_k]^2, [[a_i,a_j],a_\ell]:  1\leq j<i\leq r,\ j<k\leq r,\ k\ne i,\ \ell\in \{i,j\}\rg.
\]
By part (1) we have $K=\lg H, a_i^2 : i\in [1,r]\rg$, and it follows from   Lemma~\ref{lem:Propertyofcommutators}~(1) that each $a_i^2H$ belongs to $Z(F/H)$, so $K/H$ is abelian. Also $H/F_4\leq F'/F_4$, and  $F'/F_4$ is abelian by Lemma~\ref{lem:PropertyofF/F4}~(3), so also $H/F_4$ is abelian. 
For convenience, we let $\C_2=\{c_1,\dots,c_s\}$ and $D_K=\{d_1,\dots,d_t\}$.  Suppose that $c_1^{e_1}\cdots c_{s}^{e_{s}}d_1^{f_1}\cdots d_{t}^{f_{t}}K=K$ with $e_1,\ldots,e_s,f_1,\ldots,f_t\in\{0,1\}$. Set $g=c_1^{e_1}\cdots c_{s}^{e_{s}}d_1^{f_1}\cdots d_{t}^{f_{t}}$. Then $g\in K$, so $gH$ lies in the abelian group $K/H=\lg a_i^2H: 1\leq i\leq r\rg$, and hence $g=(a_1^2)^{\ell_1}(a_2^2)^{\ell_2}\cdots(a_{r}^2)^{\ell_{r}}h$ for some integers $\ell_1,\ell_2,\ldots, \ell_{r}$ and some $h\in H$. Also $hF_4$ lies in the abelian group $H/F_4$, and hence $h=c_1^{2i_1}\cdots c_{s}^{2i_{s}}d_1^{2j_1}\cdots d_{t}^{2j_t}h'$, where all the $i_k$ and $j_k$ are integers and $h'\in\lg F_4, [[a_i,a_j],a_\ell] : 1\leq j<i\leq r,\ \ell\in\{i,j\}\rg$.  Thus
\begin{equation*}
c_1^{e_1}\cdots c_{s}^{e_{s}}d_1^{f_1}\cdots d_{t}^{f_{t}}F_4= gF_4=(a_1^2)^{\ell_1}(a_2^2)^{\ell_2}\cdots((a_{r}^2)^{\ell_{r}})c_1^{2i_1}\cdots c_{s}^{2i_s}d_1^{2j_1}\cdots d_{t}^{2j_t}h'F_4.
\end{equation*}
However, by Theorem~\ref{Hall-basis-theorem1}, $gF_4$ has a unique representation of the form
\[
 gF_4={d'}_1^{s'_1}{d'}_2^{s'_2}\cdots {d'}_{t'}^{{s'}_{t'}} F_4,
\]
where $d'_1,d'_2,\ldots,d'_{t'}\in \cup_{u=1}^3\C_u= \{a_i,[a_i,a_j], [[a_i,a_j],a_k]: 1\leq j<i\leq r,\ j\leq k\leq r\}$ and the $s'_i$ are integers. It follows that $\ell_i=0\  (1\leq i\leq r)$, $e_k=2i_k\  (1\leq k\leq s)$, $f_k=2j_k\  (1\leq k\leq t)$ and $h'\in F_4$. Since $e_k,f_k\in\{0,1\}$, this implies that $e_1=e_2=\cdots=e_{s}=f_1=f_2=\cdots=f_t=0$. Thus $c_1K,\ldots, c_sK,d_1K\ldots, d_tK$ is a basis for 
$F'K/K\cong C_2^m$, as asserted, and $m=|\C_2\cup D_K|$. 

Finally we determine this cardinality. By Lemma~\ref{lem-basic-com1},  $|\C_2|=r(r-1)/2$, and
$|D_K|=\frac{r^3}{3} -r^2+\frac{2r}{3}$, so 
\[
m= \frac{r(r-1)}{2} + \left(\frac{r^3}{3} -r^2+\frac{2r}{3}\right) = 
\frac{r^3}{3} - \frac{r^2}{2} + \frac{r}{6} =
\frac{r(r-1)(2r-1)}{6}.
\]

This completes the proof of part (2).

(3) \quad Since $F'\leq F'K<F$, it follows that $F/F'K$ is a quotient of the abelian group $F/F'$ and hence $F/F'K$ is abelian, generated by $\{ a_iF'K : i\in [1,r]\}$. Moreover since each $a_i^2\in K\leq F'K$, the group $F/F'K\cong C_2^m$ for some $m\leq r$. Suppose that
\[
 a_1^{e_1}a_2^{e_2}\cdots a_{r}^{e_{r}}F'K=F'K,
\]
where $e_1,\ldots, e_{r}\in\{0,1\}$. Then $a_1^{e_1}a_2^{e_2}\cdots a_{r}^{e_{r}}F'=kF'$ for some $k\in K$, and it follows from part (1) that $kF'=a_1^{2f_1}a_2^{2f_2}\cdots a_{r}^{2f_{r}}F'$ for some integers $f_1,\ldots, f_{r}$. Thus $a_1^{e_1}a_2^{e_2}\cdots a_{r}^{e_{r}}F'=a_1^{2f_1}a_2^{2f_2}\cdots a_{r}^{2f_{r}}F'$. By Theorem~\ref{Hall-basis-theorem1}, $a_1F', a_2F', \ldots, a_{r}F'$ form a basis for the free abelian group $F/F'$, and this implies that $e_i=2f_i$ for each $i$. Since each $e_i=0$ or $1$, we have $e_1=e_2=\cdots=e_{r}=0$. Therefore, $a_1F'K, a_2F'K, \ldots, a_{r}F'K$ form a basis for $F/F'K$, and   $F/F'K\cong C_2^{r}$.\hfill\qed

\section{Structure of $\H(n)$}

The goal of this section is to investigate the order, the subgroups and the automorphisms of the group $\newH(n)$ in Definition~\ref{def:Hn}.
%
%
First we obtain a lower bound for the order of $\newH(n)$.

\begin{prop}\label{lem:GroupI}
Let $F, A$, and the $F_\ell$ be as in Theorem~$\ref{Hall-basis-theorem1}$, with $r=2n\geq4$, and write $A=A_0\cup B_0$, where $A_0=\{a_i : 1\leq i\leq n\}$ and $B_0=\{a_{n+i} : 1\leq i\leq n\}$. Also let $K$ be the subgroup defined in Lemma~$\ref{lem:Anormalsubgroup}$, and let $\newH(n)$ be the group defined in Definition~$\ref{def:Hn}$.  Define
\begin{align*}
 I=\lg & K, [a,a'], [b,b'],  [[a,a'],c], [[b,b'],c], [[b,a],b']\ :\ a,a'\in A_0, b,b'\in B_0, c\in A\ \rg,
\end{align*}
Then the following hold.
\begin{enumerate}
  \item [{\rm (1)}] $K < I \leq F'I=F'K$, $I\unlhd F$, and $(F/I)/(F/I)'\cong F/F'I\cong C_2^{2n}$.
  \item  [{\rm (2)}] $I$ is also generated by $K\cup D_I$, where $D_I$ is the set
  \begin{align*}
 \{ & [a_i,a_j], [a_{n+i},a_{n+j}], [[a_i,a_j],a_k], [[a_i,a_j],a_{n+\ell}], [[a_{n+i},a_{n+j}],a_{n+k}], [[a_{n+i},a_{n+j}],a_{\ell}],\\
 &[[a_{n+i'}, a_{\ell}],a_{n+j'}] : 1\leq j<i\leq n,\ j<k\leq n,\ 1\leq j'<i'\leq n,\ 1\leq\ell\leq n,\ k\neq i\}.
\end{align*}
  \item [{\rm (3)}]  $I/K\cong C_2^v$ with $v= (13n^3-15n^2+2n)/6$, and $F'I/I\cong C_2^u$ with $u=(n^3+n^2)/2$; moreover $F/I$ is nilpotent of class $3$ and order $|F/I|=2^{(n^3+n^2+4n)/2}$.
  \item [{\rm (4)}] The map  $\phi(x_i)=a_iI$ and $\phi(y_i)=a_{n+i}I$, for each $i=1,\ldots,n$, defines an epimorphism $\phi:\newH(n)\to F/I$. In particular, $|\newH(n)|\geq 2^{(n^3+n^2+4n)/2}$.
\end{enumerate}
\end{prop}

\demo (1) By definition, $K<I$, and $F'K$ contains each of the given generators for $I$, so $I\leq F'K$. This implies that $F'I\leq F'K$, and on the other hand $F'K\leq F'I$, so $F'I=F'K$. Next we prove that $I\unlhd F$. Since $F_4\leq K<I$ and $F_3/F_4\leq Z(F/F_4)$ (by Lemma~\ref{lem:PropertyofF/F4}(2)) and $K\unlhd F$ (Lemma~\ref{lem:Anormalsubgroup}), it follows that, for each $x\in F_3\cap I$, $y\in K$, and  $z\in F$, the conjugates $x^z\in xF_4\subseteq I$ and $y^z\in K<I$. Thus to prove that $I\unlhd F$,  it is sufficient to prove that, for each $a,a'\in A_0$, $b,b'\in B_0$ and $c\in A$, the conjugates $[a,a']^c, [b,b']^c$ both lie in $I$. Now $[a,a']^c = [a,a'] [[a,a'],c]$ and both factors lie in $I$, so $[a,a']^c\in I$. Similarly $[b,b']^c\in I$, and hence $I\unlhd F$. Finally $(F/I)' = F'I/I$, and hence  $(F/I)/(F/I)'\cong F/F'I = F/F'K$ which, by Lemma~\ref{lem:Anormalsubgroup}(3), is isomorphic to $C_2^{2n}$. Thus part (1) is proved.

(2) Let $I_0=\lg K, D_I\rg$, with $D_I$ as in (2). Then $I_0\leq I$ and we show that equality holds by proving that each of the given generators for $I$ lies in $I_0$. First $K$ lies in the subset so $F_4< K\leq I_0$, and it follows from Lemma~\ref{lem:Propertyofcommutators}(3) that each $[a,a']$ and $[b,b']$ lies in $I_0$. Next we consider $x:=[[a,a'],c]$. If $a=a'$ then $x=1\in I_0$, while if $c\in\{a,a'\}$, then $x\in K$ (using Lemma~\ref{lem:Anormalsubgroup} and Lemma~\ref{lem:Propertyofcommutators}(4)), and again $x\in I_0$.  Otherwise $a,a',c$ are pairwise distinct; if $c\in B_0$ then, by Lemma~\ref{lem:Propertyofcommutators}(4), $x^{\pm1}F_4$ contains an element of $D_I$ of the form $[[a_i,a_j],a_{n+\ell}]$ and hence $x\in I_0$, while if $c\in A_0$ then $x$ is of the form $x=[[a_i,a_j],a_{k}]$ with $i,j,k\leq n$ and pairwise distinct. In this case, if $\min\{i,j,k\}\ne k$ then either $x\in D_I$ or  $x^{-1}F_4$ contains $[[a_j,a_i],a_{k}]\in D_I$, by Lemma~\ref{lem:Propertyofcommutators}(4), and if $\min\{i,j,k\}= k$ then  by Lemma~\ref{lem:Propertyofcommutators}(5), $xF_4$ contains a product of two elements of $D_I$. In all cases $x\in I_0$. An analogous argument shows that each $[[b,b'],c]\in I_0$. Finally consider $x:=[[b,a],b']$. If $b'=b$ then $x\in K$ so take $b'=a_{n+j'}\ne b=a_{n+i'}$ and $a=a_\ell$. If $i'>j'$ then  $x\in D_I\subset I_0$, so we may assume that $i'< j'$. By Lemma~\ref{lem:Propertyofcommutators}(4) and (5), modulo $F_4$, $x=[[a_{n+i'},a_\ell],a_{n+j'}]$ satisfies
\[
x\equiv [[a_{n+j'},a_\ell],a_{n+i'}]\cdot [[a_{n+i'},a_{n+j'}],a_\ell]
\equiv [[a_{n+j'},a_\ell],a_{n+i'}]\cdot [[a_{n+j'},a_{n+i'}],a_\ell]^{-1},
\]
and each of these factors lies in $I_0$, so $x\in I_0$. We conclude that $I_0=I$. 

(3) Note that $F'I/I \cong F'/(F'\cap I)$, and $F'\cap I$ contains $F'\cap K$ since $F_4<K<I$. Hence $F'/(F'\cap I)$ is a quotient of $F'/(F'\cap K)$, and by Lemma~$\ref{lem:Anormalsubgroup}(2)$, $F'/(F'\cap K)\cong F'K/K\cong C_2^{r(r-1)(2r-1)/6} = C_2^{(8n^3 - 6n^2+n)/3}$ (as $|A|=r=2n$), with basis $\{c(F'\cap K): c\in \C_2\cup D_K\}$  where $D_K=\{[[a_i,a_j],a_k]: 1\leq j<i\leq 2n,\ j<k\leq 2n,\ k\ne i\}$, and $\C_2$ is as in Lemma~$\ref{lem-basic-com1}$. Since $F'K=F'I$, we deduce that $F'I/I\cong C_2^u$ and $I/K\cong C_2^v$, for some $u, v$ such that $u+v=(8n^3 - 6n^2+n)/3$. 

Let $B_2=\{[b,a]\ :\ a\in A_0, b\in B_0\}$. Then $|B_2|=n^2$, $|\C_2|=n(2n-1)$ (Lemma~\ref{lem-basic-com1}), and $\C_2$ is the disjoint union $\C_2=(\C_2\cap D_I)\cup B_2$ with $|\C_2\cap D_I|=n^2-n$.  We obtain a similar partition of $D_K$. 
Suppose that $z\in D_K$ and $z\not\in I$. Then $z=[[c,c'], c'']$ for certain $c,c',c''\in A$. First we show that $c',c''\in A_0$ and $c\in B_0$. If $c'\in B_0$, then by the definition of $D_K$ also $c,c''\in B_0$ and $z$ is an element of $D_I$ of the form $[[a_{n+i},a_{n+j}],a_{n+k}]$, which is a contradiction, so $c'\in A_0$. Next, if also $c,c''\in A_0$ then $z$ would be  an element of $D_I$ of the form $[[a_{i},a_{j}],a_{k}]$, which is a contradiction, so at least one of $c,c''$ lies in $B_0$. If $c\in A_0$ then we must have $c''\in B_0$ and then $z$ is an element of $D_I$ of the form $[[a_{i},a_{j}],a_{n+\ell}]$, which is a contradiction. Thus $c\in B_0$. If also $c''\in B_0$, then $z$ is one of the given generators for $I$ of the form $[[b,a],b']$, which is again a contradiction. Hence $c''\in A_0$, and our assertions are proved. Thus such elements $z$ have the form $z=[[a_{n+i},a_{j}],a_{k}]$, for some $i,j,k\in [1,n]$; and as $z\in D_K$ we have $j<k$. Let
\[
B_3:=\{ [[a_{n+i},a_{j}],a_{k}]\ :\ 1\leq j<k\leq n,\ 1\leq i\leq n \}.
\]
Then by Lemma~\ref{lem:arith}, $|B_3|=n\cdot n(n-1)/2 = (n^3-n^2)/2$, and $D_K$ is the disjoint union $D_K=B_3'\cup B_3$ with $B_3'\subset I$. Now  $|D_K|=\frac{8n^3}{3} -4n^2+\frac{4n}{3}$ (Lemma~\ref{lem:Anormalsubgroup}), so $|B_3'|=|D_K|-|B_3|= (13n^3-21n^2+8n)/6$.

An analogous argument to that in the proof of Lemma~\ref{lem:Anormalsubgroup}(2) shows that $\{cI\ :\ c\in B_2\cup B_3\}$ forms a basis for $F'I/I\cong C_2^u$ and hence $u=|B_2|+|B_3| = n^2 + (n^3-n^2)/2 = (n^3+n^2)/2$. Thus by part (1), $|F/I|=2^{2n+u}$ and $2n+u = (n^3+n^2+4n)/2$. Also $I/K\cong C_2^v$ with $v= |\C_2\cap I|+|B_3'|=(n^2-n)+(13n^3-21n^2+8n)/6 = (13n^3-15n^2+2n)/6$. 

To see that $F/I$ is nilpotent of class 3, we note first that its derived subgroup is $(F/I)'=F'I/I\cong C_2^{2n}$. The second term in the lower central series is $F_3I/I$. It follows from Lemma~\ref{lem:PropertyofF/F4}(1) that $F'I=\lg \C_2, F_3, I\rg$, and hence $F'I/F_3I$ is generated by $\{ cF_3I : c\in\C_2, c\not\in I\}$. We showed above that $\C_2\setminus I = B_2$ has size $n^2$, and hence  $F'I/F_3I=C_2^{u_2}$ with $u_2\leq n^2<u$. The third term in the lower central series is $F_4I/I$, which is trivial since $F_4<I$. Thus $F/I$ is nilpotent of class 3, as asserted. This completes the proof of part (3).

(4)  Setting $\phi(x_i)=a_iI$ and $\phi(y_i)=a_{n+i}I$, for each $i=1,\ldots,n$, we have a map from the generating set for $\newH(n)$ in Definition~\ref{def:Hn} to the set $\{aI\ :\ a\in A\}$ of generators for $F/I$. Moreover, extending $\phi$ to a map on words in these generators of $\newH(n)$, for each relator  
$w\in \mathcal{R}$ in Definition~\ref{def:Hn}, $\phi(w) = w'I$ such that either $w'\in K$ or $w'$ is one of the given generators for $I$. Thus the images $\phi(x_i), \phi(y_i)$ $(1\leq i\leq n)$ satisfy all the given relations of $\newH(n)$, and hence, by von Dyck's Theorem (see \cite[Theorem~2.2.1]{Robinson}, the extension of $\phi$ to $\newH(n)\to F/I$ is an epimorphism. This  completes the proof of the proposition.
\hfill\qed 

Our next task is to prove that the epimorphism $\phi$ in Proposition~\ref{lem:GroupI}(4) is in fact an isomorphism. We need the following information about certain commutators in $\newH(n)$. 

\begin{lem}\label{lem:H-props}
    Let $\newH(n) = \lg X_0\cup Y_0 \mid \mathcal{R}\rg$ be the group defined in Definition~$\ref{def:Hn}$, where $n\geq2$. 
Then, 
\begin{enumerate}
    \item [{\rm(1)}] for all $z,z'\in X_0\cup Y_0$, $[[z,z'],z]=[[z,z'],z']=1$;
    \item [{\rm(2)}] for $z,z',z''\in X_0\cup Y_0$,  we have $[[z,z'],z'']\in Z(\newH(n))$, $[z,z']=[z',z]$, and  $[z,z']^2=[[z,z'],z'']^2=1$. 
    \item [{\rm(3)}] the third term $\newH(n)_4$ of the lower central series for $\newH(n)$ is trivial, so $\newH(n)$ is nilpotent of class at most $3$.
\end{enumerate}
\end{lem}

\demo
(1)  We use the first few relations in $\mathcal{R}$. Firstly, $z^2=(z')^2=1$, so $(zz')^2=[z,z']$. If both $z,z'$ lie in $X_0$ or both lie in $Y_0$, then we have the relation $[z,z']=1$, and hence $[[z,z'],z]=[[z,z'],z']=1$. 
So we may assume that $z\in X_0$, say, and $z'\in Y_0$. Then we have $(zz')^4= [z,z']^2=1$, and hence $\lg z,z'\rg=D_8$ or $C_2^2$. In either case $[z,z']=(zz')^2$ is centralised by $z$ and $z'$, and this implies that $[[z,z'],z]=[[z,z'],z']=1$. This proves part (1). 

(2) Let $u=[[z,z'],z'']$. Then for all $z'''\in X_0\cup Y_0$, $[u,z''']\in\mathcal{R}$ and hence $[u,z''']=1$.
Since $[u,z''']=u^{-1} u^{z'''}$, this implies that $u=u^{z'''}$, and hence $u\in Z(\newH(n))$. If both $z,z'$ lie in $X_0$ or both lie in $Y_0$, then $[z,z']=[z',z]=1$  and the other assertions follow. If $z\in X_0$ and $z'\in Y_0$, then $[z,z']^2\in\mathcal{R}$ and hence $[z,z']^2=1$, and this implies that $[z,z']^{-1}=[z,z']$. However $[z,z']^{-1}=[z',z]$, and hence $[z,z']=[z',z]$. Also $u^2\in\mathcal{R}$ in this case and so $u^2=1$. Finally, if $z\in Y_0$ and $z'\in X_0$, then $[z',z]^2\in\mathcal{R}$, and the argument just given shows that  $[z,z']=[z',z]$ and $[z,z']^2=1$. This implies that $u^2=[[z,z'],z'']^2=[[z',z],z'']^2$, which lies in $\mathcal{R}$ and hence is trivial. 

(3) By part (2), for all $z,z',z'', z'''\in X_0\cup Y_0$,  we have $[[z,z'],z'']\in Z(\newH(n))$ and hence $[[[z,z'],z''],z'''] =1$. On the other hand, by \cite[Hilfsatz 1.11(a)]{Huppert}, $\newH(n)_4$ is generated by the set of all conjugates $[[[z,z'],z''],z''']^h$ of such commutators  by elements  $h\in\newH(n)$, and hence $\newH(n)_4$ is trivial. Thus $\newH(n)$ is nilpotent of class at most $3$.
\hfill\qed


\begin{prop}\label{lem:Isomorphism}
   Let $\newH(n) = \lg X_0\cup Y_0 \mid \mathcal{R}\rg$ be the group defined in Definition~$\ref{def:Hn}$, where $n\geq2$, and let $F, A=A_0\cup B_0, I, K$ be as in Proposition~$\ref{lem:GroupI}$. Then the map $\psi: AI/I\to X_0\cup Y_0$ such that  $\psi(a_iI)=x_i$ and $\psi(a_{n+i}I)=y_i$, for $i=1,\dots,n$, defines an epimorphism $\psi:F/I\to \newH(n)$, such that $\psi$ is the inverse of the map $\phi$ in Proposition~$\ref{lem:GroupI}(4)$. In particular  $\newH(n)\cong F/I$ and $|\newH(n)|= 2^{(n^3+n^2+4n)/2}$.
\end{prop}

\begin{rem}
It follows from Proposition~$\ref{lem:GroupI}(1)$ that $F/I$ is isomorphic to the group with presentation $\olF=\lg A_0\cup B_0 \mid F_4\cup B_K\cup D_I\rg$, where $B_K, D_I$ are as in Lemma~$\ref{lem:Anormalsubgroup}$ and Proposition~$\ref{lem:GroupI}(2)$. In the proof we will work with the group $\olF$.
\end{rem}

\demo
 Interpreting $F/I$ as the group $\olF:=\lg A_0\cup B_0 \mid F_4\cup B_K\cup D_I\rg$, the map  $\psi$ becomes  a bijection $\psi:A_0\cap B_0\to X_0\cup Y_0$ given by $\psi: a_i\to x_i, a_{n+i}\to y_i$, for $i=1,\dots,n$. Consider the extension of $\psi$ to a map on words in these generators of $\olF$, so that for each element (relator) $w\in F_4\cap B_K\cup D_I$, $\psi(w)$ is the same word in $X_0\cup Y_0$. We check that $\psi(w)$ is equal to the identity of $\newH(n)$ for each $w$ and then apply von Dyck's Theorem \cite[Theorem 2.2.1]{Robinson}. If $w\in F_4$ then $\psi(w)=1$ since $\newH(n)$ has nilpotency class at most $3$ by Lemma~\ref{lem:H-props}(2). 
 
 Next we consider the elements $w$ of $B_K$ as in Lemma~\ref{lem:Anormalsubgroup}(1). If $w=a^2$ for $a\in A_0\cup B_0$, then $\psi(w)\in\mathcal{R}$ and so $\psi(w)=1$. If $w=[a,a']^2$ or $w=[[a,a'],a'']^2$, for $a,a', a''\in A_0\cup B_0$, then  $\psi(w)=1$ by Lemma~\ref{lem:H-props}(3). 
 Finally suppose that $w=[[a,a'],a'']$, for $a,a'\in A_0\cup B_0$ with $a''\in\{a,a'\}$. Then it follows from Lemma~\ref{lem:H-props}(1), that $\psi(w)=1$. 
 
 Finally we consider the elements $w$ of $D_I$ as in Proposition~$\ref{lem:GroupI}(2)$. If $w=[a,a']$ for $a, a'\in A_0$ or $a,a'\in B_0$, then $\psi(w)\in\mathcal{R}$ and so $\psi(w)=1$. If $w=[[a,a'],a'']$, for $a,a', a''\in A_0\cup B_0$ such that either $a,a'\in A_0$ or $a,a'\in B_0$, then again $\psi(w)=1$ since in these cases $\psi([a,a'])\in\mathcal{R}$. Finally we consider elements $w=[[a,a'],a'']$, for $a,a''\in B_0$ and $a'\in A_0$. For these elements, $\psi(w)\in\mathcal{R}$ and so $\psi(w)=1$. 

 It now follows from von Dyck's Theorem \cite[Theorem 2.2.1]{Robinson} that this map $\psi$ defines an epimorphism $\olF\to \newH(n)$.  By the definitions,  $\psi$ is the inverse of the map $\phi$ in Proposition~$\ref{lem:GroupI}(4)$, and hence $\newH(n)\cong F/I$. The order $|\newH(n)| = |F/I|$ follows from Proposition~\ref{lem:GroupI}(3).
\hfill\qed

\subsection{Subgroups and automorphisms of $\newH(n)$}

The main purpose of this subsection is to prove the following theorem.

\begin{thm}\label{prop-2groupNew}
Let $\newH(n)$ be as in Definition~\ref{def:Hn} and let $\newH(n)_3$ be the third term in the lower centre series of $\newH(n)$. Then the following hold.
\begin{enumerate}
\item [{\rm(1)}] For arbitrary $x,x'\in X_0$, $y\in Y_0$, we have $[[x,y],x']=[[x',y],x]$. 
  \item [{\rm(2)}] $\newH(n)/(\newH(n)')\cong C_2^{2n}$, and $\newH(n)$ is an $n$-dimensional mixed 
  dihedral group relative to $X :=\lg X_0\rg$ and $Y :=\lg Y_0\rg$.
    \item [{\rm(3)}]   $\newH(n)'=W\times \newH(n)_3\cong C_2^{n^2(n+1)/2}$, where 
    \begin{enumerate}
        \item [{\rm(a)}] $W=\lg [x_i,y_j]: 1\leq i,j\leq n\rg \cong C_2^{n^2}$, and 
        \item [{\rm(b)}] $\newH(n)_3=\lg [[x_i,y_j],x_k]: 1\leq i,j\leq n,\ i<k\leq n\rg \cong C_2^{n^2(n-1)/2}$; and moreover $\newH(n)_3\leq Z(\newH(n))$. 
    \end{enumerate}
    %
  \item [{\rm(4)}] For any $a\in \newH(n)$ and $b,b'\in Y$ we have $[[b,a],b']=1$ and $[[a,b],b']=1$.
  \item [{\rm(5)}] For any $g\in \Aut(X)\times\Aut(Y)$, $g$ induces an automorphism of $\newH(n)$.
  \item [{\rm(6)}] Let $c\in X$ and let $d\in Y$ so that $c=x_{i_1}x_{i_2}\dots x_{i_k}$ and $d=y_{j_1}y_{j_2}\dots y_{j_\ell}$ for some $x_{i_1},x_{i_2},\dots, x_{i_k}\in X_0$ and $y_{j_1}, y_{j_2},\dots, y_{j_\ell}\in Y_0$, with these expressions chosen so that $k, \ell$ are minimal. Then
  \[
   [c,d]\newH(n)_3= \prod_{1\leq u\leq k, 1\leq v\leq \ell} [x_{i_u},y_{j_v}] \newH(n)_3.
  \]
\end{enumerate}
\end{thm}

\f\demo 
(1) Let $x,x'\in X_0$ and $y\in Y_0$. By Lemma~\ref{lem:H-props} (2), any weight 3 commutator involving $x,x',y$ is in the center of $\newH(n)$. Thus, by Lemma~\ref{Witt-formula}, we have
\[
[[x,y^{-1}],x']\cdot[[y,(x')^{-1}],x]\cdot [[x',x^{-1}],y]=1.
\]
Since $x^2=(x')^2=y^2=1$ and $[x,x']=1$ the above equation becomes $[[x,y],x']\cdot [[y,x'],x]=1$, and hence $[[x,y],x']=[[y,x'],x]^{-1}$. Using the relations $[[x',y],x]^2=1$ and $[x',y]^2=1$ (which implies that $[x',y]=[y,x']$) we have $[[y,x'],x]^{-1}=[[x',y],x]^{-1}=[[x',y],x]$, and part (1) holds. 

(2) It follows from Proposition~\ref{lem:Isomorphism} that $\newH(n)\cong F/I$, and from Proposition~\ref{lem:GroupI} (1) that $\newH(n)/\newH(n)'\cong F/F'I\cong C_2^{2n}$. Also, $\newH(n)=\lg X,Y\rg$, $X 
=\lg X_0\rg\cong C_2^n$ and $Y =\lg Y_0\rg\cong C_2^n$, by Definition~\ref{def:Hn}, and hence, by Definition~\ref{mix-dih}(a), $\newH(n)$ is an $n$-dimensional mixed dihedral group relative to $X$ and $Y$.

(3) By Propositions~\ref{lem:Isomorphism} and~\ref{lem:GroupI} (3),  
 $\newH(n)'\cong F'I/I \cong C_2^u$ with $u=(n^3+n^2)/2$, and by \cite[Hilfsatz 1.11(a) and (b)]{Huppert}, 
\[
 \newH(n)'=\lg [z,z'],[[z,z'],z'']^h\mid z,z',z''\in X_0\cup Y_0,\ h\in\newH(n)\rg.
\]
Also, by Lemma~\ref{lem:H-props} (2), $\newH(n)_3\leq Z(\newH(n))$ and hence $[[z,z'],z'']^h=[[z,z'],z'']$ for each $h\in\newH(n)$. Let $x,x'\in X_0$ and $y,y'\in Y_0$. Since $[x,x']= [y,y']= 1$ are relations in $\R$, the only weight two generators required are $[x_i,y_j]$, for $1\leq i,j\leq n$, and for the weight three generators $[[z,z'],z'']$, we may assume that one of $z,z'$ lies in $X_0$ and the other lies in $Y_0$. 
Since $\newH(n)_4=1$ by Lemma~\ref{lem:H-props}(3), it follows from Lemma~\ref{lem:Propertyofcommutators}(4) that $[[x,y],y']^{-1}=[[y,x],y']$, and since $[[y,x],y']=1$ is a relation in $\R$, also $[[x,y],y']=1$, and so  the only weight three generators required are those of the form $[[y,x],x']$ or $[[x,y],x']$.  
Further by part (1), $[[x,y],x']=[[x',y],x]$, and by Lemma~\ref{lem:H-props} (1), $[[x,y],x]=1$. Thus we have
\[
 \newH(n)'=\lg [x_i,y_j],[[x_i,y_j],x_k]\mid 1\leq i,j\leq n,\ i<k\leq n\rg.
\]
Since there are precisely $u=(n^3+n^2)/2$ generators in the generating set above, and since $\newH(n)'=C_2^u$, we conclude that $\newH(n)'=W\times \newH(n)_3$, where $W=\lg [x_i,y_j]: 1\leq i,j\leq n\rg\cong C_2^{n^2}$ and $\newH(n)_3=\lg [[x_i,y_j],x_k]: 1\leq i,j\leq n,\ i<k\leq n\rg\cong C_2^{n^2(n-1)/2}$. 


(4) By part (3), we have $[b,a]=gh$ where $g \in W$ and $h \in \newH(n)_3$. Moreover, by part~(3) we also have $\newH(n)_3\leq Z(\newH(n))$. This implies that 
\[
[[b,a],b']=[gh,b']=(gh)^{-1}(b')^{-1}(gh)b'=g^{-1}(b')^{-1}gb'=[g,b'].
\]
Again by part (3), we have $g=w_1w_2\dots w_s$ for some $w_1,w_2,\dots,w_s\in \{[x_i,y_j]: 1\leq i,j\leq n\}$. We will use induction on $s$ to prove that $[g,y_k]=1$, for any $y_k\in Y_0$. In the proof of part (3) we showed that $[[x_i,y_j],y_k]=1$, for all $i,j,k$. 
Thus, if $s=1$, then $[g,y_k]=[w_1,y_k]=1$. Now assume that $s>1$, and assume inductively that $[g,y_k]=1$ if $g\in W$ can be expressed as a word of length less than $s$ in the generators. Then, by \cite[Hilfssatz III.1.2(c)]{Huppert},
\[
 [g,y_k]=[w_1w_2\dots w_{s-1}w_s,y_k]=[w_1w_2\dots w_{s-1},y_k]^{w_s}[w_s,y_k],
\]
and also $[w_s,y_k]=1$ from the case $s=1$. 
Since by part (3) the commutator subgroup $\newH(n)'$ is abelian, we have
\[
[g,y_k]= [w_1w_2\dots w_{s-1},y_k]^{w_s}=[w_1w_2\dots w_{s-1},y_k].
\]
By induction $[w_1w_2\dots w_{s-1},y_k]=1$, and hence $[g,y_k]=1$. This implies that $g$ commutes with every $y_k \in Y_0$. Since $Y$ is generated by $Y_0$, it follows that $g$ centralises $Y$, and hence we have $[g,b']=1$. Thus, $[[b,a],b']=1$. 
By part~(3), $\newH(n)'$ is an elementary abelian $2$-group, so $[b,a]=[b,a]^{-1} = [a,b]$, and hence also $[[a,b],b']=1$.

(5) Let $(g_1,g_2)\in \Aut(X)\times\Aut(Y)$, and let $x_i'=x_i^{g_1}$ and $y_i'=y_i^{g_2}$, for each $i\in\{1,\ldots,n\}$. Set $X_0'=\{x_i': i=1,2,\ldots,n\}$ and $Y_0'=\{y_i': i=1,2,\ldots,n\}$. We will apply von Dyck's Theorem (see \cite[Theorem~2.2.1]{Robinson}) to show that the map $\phi:X_0\cap Y_0\to \newH(n)$ given by $x_i\to x_i',\ y_i\to y_i'$, for $i=1,\dots,n$, extends uniquely to an automorphism $\phi$ of $\H(n)$.
To do this, it is sufficient to show that $\newH(n)$ is generated by $X_0'\cup Y_0'$ and that for every relation $w(x_1,\ldots,x_n,y_1,\ldots,y_n)=1$ in $\mathcal{R}$, we have $w(x_1',\ldots,x_n',y_1',\ldots,y_n')=1$. First, $X=\lg X_0\rg\cong Y=\lg Y_0\rg\cong C_2^{n}$ and $\newH(n)=\lg X,Y\rg$, by the definition of $\newH(n)$. Also $\lg X_0'\rg=X$ and $\lg Y_0'\rg=Y$, by the definition of $g_1$ and $g_2$. Hence $\newH(n)=\lg X'_0\cup Y'_0\rg$. Next we consider the relations. Let $a,a'\in X_0', b,b'\in Y_0'$ and $c, c'\in X_0'\cup Y_0'$. Then $c^2=1$ and $[a,a']=[b,b']=1$. By part (3), we have $[a,b],[[a,b],c]\in\newH(n)'\cong C_2^{n^2(n+1)/2}$, and thus, $[a,b]^2=1$ and $[[a,b],c]^2=1$. By part (4), $[[b,a],b']=1$. Finally, $[[[a,b],c],c']$ lies in $\newH(n)_4$ and hence is trivial, by Lemma~\ref{lem:H-props}(3). This proves part~(5).

(6) In the following we write, for convenience, $w\equiv w'\pmod{\newH(n)_3}$ if and only if $w\newH(n)_3 = w'\newH(n)_3$. 
First we apply \cite[Hilfssatz III.1.2(c)]{Huppert} several times: $[gh,f]=[g,f]^h\cdot 
[h,f]=[g,f]\cdot [[g,f],h]\cdot [h,f]$, for any $g,h,f\in\newH(n)$. Writing $c=c'x_{i_k}$, this implies that $[c,d]=[c',d]\cdot [[c',d],x_{i_k}]\cdot [x_{i_k},d]$, and since $[[c',d],x_{i_k}]\in Z(\newH(n)_3)$ (by part~(3)), it follows that $[c,d]\equiv [c',d]\cdot  [x_{i_k},d] \pmod{\newH(n)_3}$. Repeating this $k$ times we obtain
\[
[c,d]\equiv [x_{i_1},d]\dots  [x_{i_k},d] \pmod{\newH(n)_3}.
\]
Now we apply \cite[Hilfssatz III.1.2(b)]{Huppert} several times: $[g,hf]=[g,f]\cdot[g,h]^f=[g,f]\cdot [g,h]\cdot [[g,h],f]$, for any $g,h,f\in\newH(n)$. Writing $d=d'y_{j_\ell}$, this implies, for all $u$, that $[x_{i_u},d]=[x_{i_u},y_{j_\ell}]\cdot [x_{i_u},d']\cdot [[x_{i_u},d'],y_{j_\ell}]$, and since $[[x_{i_u},d'],y_{j_\ell}]=1$ (by part~(4)), it follows that $[x_{i_u},d]=[x_{i_u},y_{j_\ell}]\cdot [x_{i_u},d']$. Repeating this $\ell$ times for each $u$, and using the fact that $\newH(n)'/\newH(n)_3$ is abelian (by part~(3)), we obtain
\[
[c,d]\equiv \prod_{1\leq u\leq k, 1\leq v\leq \ell} [x_{i_u},y_{j_v}] \pmod{\newH(n)_3}.
\]
This completes the proof.
\hfill\qed

\section{Proof of Theorem~\ref{t:semi-sym}}

Let $\newH(n)$, $X=\lg X_0\rg$, $Y=\lg Y_0\rg$ and $\R$ be as in Definition~\ref{def:Hn}. By Theorem~\ref{prop-2groupNew}~(2), $\newH(n)$ is an $n$-dimensional mixed dihedral group relative to $X$ and $Y$. By Proposition~\ref{lem:Isomorphism}, the order of $\newH(n)$ is $2^{(n^3+n^2+4n)/2}$. Let $\Gamma=C(\newH(n),X,Y)$ and $\Sigma=\Sigma(\newH(n),X,Y)$, as in Definition~\ref{mix-dih}. It follows from Lemma~\ref{lem:prop-mixed-dih}(5) that $\Sigma$ has valency $2^n$, and from Definition~\ref{mix-dih}(b) that $|V(\Sigma|=2\cdot|\newH(n):X| = 2^a$, where $a=1 + (n^3+n^2+4n)/2 - n = 1+ (n^3+n^2+2n)/2$. It remains  for us to prove that $\Sigma$ is semisymmetric and locally $2$-arc-transitive.

First we prove that $\Sigma$ is locally $2$-arc-transitive. By Lemma~\ref{lem:prop-mixed-dih}(1), (3) and (4), $\Sigma$ is the clique graph of $\Gamma$, $\Aut(\Gamma)=\Aut(\Sigma)$ contains $G:= \newH(n) \rtimes A(\newH(n),X,Y)$,  the group $\newH(n)$ has two obits on $V(\Sigma)$, namely
$\{Xh : h\in \newH(n)\}$ and  $\{Yh : h\in \newH(n)\}$, and $\newH(n)$ acts regularly on $E(\Sigma)$. Further the stabiliser in $G$ of the $1$-arc $(X,Y)$ of $\Sigma$ is the subgroup $A(\newH(n),X,Y)$. By Theorem~\ref{prop-2groupNew}(5), $A(\newH(n),X,Y)$ contains $\Aut(X)\times\Aut(Y)$.  By \eqref{eq-2}, $(X,Y,Z)$ is a $2$-arc of $\Sigma$ if and only if $Z=Xz$ for some $z\in \newH(n)$ such that $Xz\cap Y\ne \emptyset$. Thus $Z=Xy$ for some $y\in Y$ and since $Z\ne X$, we have $y\ne 1$. Since $\Aut(Y)\cong \GL_n(2)$ is transitive on $Y\setminus\{1\}$ it follows that $\Aut(X)\times\Aut(Y)$ is transitive on all the $2$-arcs of the form $(X,Y,Z)$, and hence the stabiliser in $G$ of $X$ is transitive on all the $2$-arcs of $\Sigma$ with first vertex $X$. An analogous argument with $X$ and $Y$ interchanged shows that the stabiliser in $G$ of $Y$ is transitive on all the $2$-arcs of $\Sigma$ with first vertex $Y$, and it follows that $\Sigma$ is locally $2$-arc-transitive.


Showing that $\Sigma$ is semisymmetric is the most delicate part of the proof. In the smallest case, where $n=2$, a computation using Magma~\cite{BCP} shows that $\Sigma$ is semisymmetric (see Remark~\ref{rem-magma} for a description of these computations). Thus we assume that $n\geq3$. By Lemma~\ref{lem:prop-mixed-dih} (4), $\Sigma$ is edge-transitive. Thus, to show that $\Sigma$ is semisymmetric it is sufficient to prove that $\Aut(\Sigma)$ is not transitive on $V(\Sigma)$. We suppose to the contrary that $\Aut(\Sigma)$ is transitive on $V(\Sigma)$, and seek a contradiction. Under this assumption $\Sigma$ is a $2$-arc-transitive graph of order a $2$-power and valency $2^n\geq8$. We shall process the proof by the following four steps.\medskip

\f{\bf Step~1.}\ $\newH(n)'\unlhd \Aut(\Sigma)$.\smallskip

Let $u=X\in V(\Sigma)$, and let $A:= \newH(n) \rtimes (\Aut(X)\times\Aut(Y))$. Then by \eqref{eq-2}, $\Sigma(u) = \{Yx : x\in X\}$, and hence by Lemma~\ref{lem:prop-mixed-dih}(4) and Theorem~\ref{prop-2groupNew}(5), the kernel of the action of $\Aut(\sigma)_u$ on $\Sigma(u)$ contains $\Aut(Y)$. Thus Lemma~\ref{reducetoK2n2n} applies, and so there exists a $2$-group $M\unlhd \Aut(\Sigma)$ such that $M\leq \Aut(\Sigma)^+$, $M$ is semiregular on $V(\Sigma)$, and $\Sigma$ is an $M$-normal cover of $\Sigma_M\cong \K_{2^n,2^n}$.  As noted above  $\newH(n)\unlhd A\leq  \Aut(\Sigma)^+$ and $\newH(n)$ acts regularly on $E(\Sigma)$ and $A$ is locally $2$-arc-transitive on $\Sigma$ (by Lemma~\ref{lem:prop-mixed-dih}(4)). Hence $M\newH(n)\unlhd MA\leq \Aut(\Sigma)^+$, $M\newH(n)$ is a $2$-group (since both $M$  and $\newH(n)$ are $2$-groups), and $M\newH(n)$ is edge-transitive on $\Sigma$ (since $\newH(n)$ is transitive on $E(\Sigma)$), and its vertex-orbits are the two biparts of $\Sigma$. Let $\Phi$ be the Frattini subgroup of $M\newH(n)$, so $\Phi$ is a characteristic subgroup of $M\newH(n)$ and hence $\Phi\unlhd MA$. 
If $\Phi$ were transitive on one of the biparts, say $O$, of $\Sigma_M$, and if $v\in O$, then $M\newH(n) = (M\newH(n))_v\Phi$, and by the properties of a Frattini subgroup (\cite[Satz III.3.2(a)]{Huppert}),  $M\newH(n) = (M\newH(n))_v$, contradicting the    fact that $M\newH(n)$ is transitive on each bipart of $\Sigma$. Thus $\Phi$ is intransitive on each bipart of $\Sigma_M$, and $\Phi\unlhd MA$. On the other hand
$\Sigma$ is locally $(MA,2)$-arc transitive (since it is locally $(A,2)$-arc transitive), and hence 
$\Sigma_M$ is locally $(MA/M,2)$-arc transitive (by \cite[Lemma 5.1]{GLP-tamc}). 

Since $\Sigma_M\cong \K_{2^n,2^n}$, this means that $MA$ acts $2$-transitively on each bipart $O$ of $\Sigma_M$, and since $\Phi$ is an intransitive normal subgroup of $MA$ it follows that $\Phi$ acts trivially on $O$, for each bipart $O$  of $\Sigma_M$. Hence $\Phi$ is contained in the kernel of the action of $MA$ on $\Sigma_M$, that is, $\Phi\leq M$.  Thus $M\newH(n)/M$ is a quotient of $M\newH(n)/\Phi$ and hence $M\newH(n)/M\cong C_2^s$ for some $s$ (by Lemma~\ref{basis}). Since $M\newH(n)$ is edge-transitive on $\Sigma$, it follows that $M\newH(n)/M$ is an abelian group acting transitively  on $E(\Sigma_M)$, and hence $M\newH(n)/M$ is regular  on $E(\Sigma_M)$ (see \cite[Lemma 2.4]{PS}), so $s=2n$ and $\newH(n)/(M\cap \newH(n))\cong C_2^{2n}$. The group induced by $A$ on $\Sigma_M$ is $A/(A\cap M)\cong MA/M$,  and is isomorphic to $C_2^{2n} \rtimes (\Aut(X)\times\Aut(Y))$.
Thus both $A$ and $MA$ are edge-transitive on $\Sigma$ with edge-stabilisers isomorphic to $\Aut(X)\times\Aut(Y)$, and hence $MA=A$, that is, $M\leq A$. Then since $\newH(n)$ is the largest normal $2$-subgroup of $A$, we have $M\leq \newH(n)$. 
It follows that $\newH(n)/M\cong C_2^{2n}$ and hence $M\leq \newH(n)'$; and since $\newH(n)/\newH(n)'\cong C_2^{2n}$ by Theorem~\ref{prop-2groupNew}(2), we conclude that $M=\newH(n)'$, and the assertion of Step 1 is proved.

\medskip For the next part of the argument we exploit the fact that $\Aut(\Sigma)=\Aut(\Gamma)$, recalling that $\Gamma=C(\newH(n),X,Y)$ is the Cayley graph $\Cay(\newH(n), S)$, where $S=(X\cup Y)\setminus\{1\}$. 
We will frequently use the basic fact \cite[Lemma 2.1]{HPZ} about mixed dihedral groups  that the natural projection map $\phi:h\to h\newH(n)'$ determines an isomorphism $\newH(n)/\newH(n)'\cong \phi(X)\times\phi(Y)\cong X \times Y$.  In Step~2 we study the subset $\Gamma_4(1)$ of the vertex set $\newH(n)$ of $\Gamma$ consising of all elements $h$ which can be reached by a path of length at most four from the vertex $1$, so  $\Gamma_4(1)$ consists of all elements $h\in\newH(n)$ such that $h=h_1h_2\dots h_k$ with each $h_i\in S$ and $0\leq k\leq 4$.\smallskip

\medskip
\f{\bf Step~2.}\ $\G_4(1)\cap \newH(n)'=\{ 1\}\cup S'$, where $S':= \{[x,y]: x\in X\setminus\{1\}, y\in Y\setminus\{1\}\}.$\smallskip

Let $\G_4(1)':=\G_4(1)\cap \newH(n)'$. Clearly $1\in \G_4(1)'$ and $\G_4(1)'\subseteq \newH(n)'$. Suppose that $h\in \G_4(1)'\setminus\{1\}$, so $h=h_1h_2\dots h_k$ with each $h_i\in S$ and $1\leq k\leq 4$. Choose such an expression for $h$ with $k$ minimal. Note in particular that $h\in\newH(n)'$ and hence $\phi(h)=1$, with $\phi$ as above. If all the $h_i\in X\setminus\{1\}$ then $h\in X$, and since $h\ne1$ it follows from \cite[Lemma 2.1]{HPZ} that $\phi(h)\ne1$ which is a contradiction. We obtain a similar contradiction if all the $h_i$ lie in $Y\setminus\{1\}$. Hence $2\leq k\leq 4$ and not all the $h_i$ lie in the same set, $X\setminus\{1\}$ or $Y\setminus\{1\}$. Next if there exists a unique $i$ such that $h_i\in X\setminus\{1\}$, then $\phi(h)=\phi(h_i)\cdot a$ for some $a\in\phi(Y)$, and again we find that $\phi(h)\ne 1$, and obtain a contradiction. Thus at least two of the $h_i$ lie in $X\setminus\{1\}$ and, similarly, at least two of the $h_i$ lie in $Y\setminus\{1\}$. This means that $k=4$, and exactly two of the $h_i$ lie in $X\setminus\{1\}$, say $x$ and $x'$, and exactly two of the $h_i$ lie in $Y\setminus\{1\}$, say $y$ and $y'$. Then $\phi(h)=\phi(xx')\cdot \phi(yy')$. If at least one of $xx'$ or $yy'$ is nontrivial then $\phi(h)\ne 1$ by \cite[Lemma 2.1]{HPZ}, and we have a contradiction. Thus $x'=x^{-1}=x$ and $y'=y^{-1}=y$.  Further, if $h_i=h_{i+1}$ for some $i$ we would have $h_ih_{i+1}=1$ and obtain a shorter expression for $h$. Thus the minimality of $k$ implies that $h=xyxy=[x,y]$, or $h=yxyx=[y,x]=[x,y]$ (where the last equality uses the facts that each of $x, y$ and $[x,y]$ is equal to its inverse). Thus Step 2 is proved. 
\medskip

\medskip
\f{\bf Step~3.}\  $\Aut(\Sigma)_1$ fixes setwise the subset  $S'$ of $V(\G)$ in Step 2, and acts transitively on $S'$.
\smallskip

By Step~1, we have $\newH(n)'\unlhd \Aut(\Sigma)$. Thus $\a^{-1}h\a\in \newH(n)'$ for all $\a$ in the vertex stabiliser $\Aut(\Sigma)_1$ and $h\in \newH(n)'$. Since $\newH(n)$ acts on $V(\G)=\newH(n)$ by right multiplication, the image of the vertex $h\in\newH(n)'$ under $\a\in\Aut(\Sigma)_1$ is
\[
h^\a=(1^h)^\a=1^{h\a}=1^{\a^{-1}h\a}=\a^{-1}h\a.
\] 
This implies that $\Aut(\Sigma)_1$ fixes setwise  the subset $\newH(n)'$ of $V(\G)$.
Since $\Aut(\Sigma)_1$ also  fixes $\G_4(1)$ setwise, and fixes the vertex $1$, it follows that $\Aut(\Sigma)_1$ fixes $(\G_4(1)\cap \newH(n)')\setminus\{1\}$ setwise.
By Step 2, we have $(\G_4(1)\cap \newH(n)')\setminus\{1\} = \{[x,y]: x\in X\setminus\{1\}, y\in Y\setminus\{1\}\}=S'$. Recall that $\Aut(X)\times \Aut(Y)\leq \Aut(\Sigma)_1$ and, since $\Aut(X)\times \Aut(Y)$ normalises $\newH(n)$, that $\Aut(X)\times \Aut(Y)$ acts on $V(\G)=\newH(n)$ via its natural action. Since $\Aut(X)$ is transitive on $X\setminus\{1\}$ and $\Aut(Y)$ is transitive on $Y\setminus\{1\}$, it follows that $\Aut(X)\times \Aut(Y)$, and hence also $\Aut(\Sigma)_1$, is transitive on $S'$. Thus Step~3 is proved.

\medskip
\f{\bf Step~4.}\ A final contradiction.
\smallskip

For the final part of the proof we analyse a Cayley graph related to $\Gamma$, namely the graph $\Lambda :=\Cay(\newH(n), S\cup S')$. Note that the right multiplication action of $\newH(n)$ yields $\newH(n)$ as a subgroup of automorphisms of the graphs $\G :=\Cay(\newH(n), S)$ and $\Cay(\newH(n), S')$, and hence also $\newH(n)\leq \Aut(\Lambda)$. Moreover, since $\Aut(\Sigma)=\Aut(\G)$, the group $\Aut(\Sigma)_1$, in its natural action on $V(\Lambda)=\newH(n)$, leaves $S$ invariant, and by Step~3, $\Aut(\Sigma)_1$ also leaves $S'$ invariant  (and is transitive on it), and hence $\Aut(\Sigma)_1$ leaves $S\cup S'$ invariant. Therefore also $\Aut(\Sigma)_1\leq \Aut(\Lambda)$ and hence, since $\Aut(\Sigma)=\newH(n)\Aut(\Sigma)_1$, we have $\Aut(\Sigma)\leq \Aut(\Lambda)$. Now $\Lambda(1)=S\cup S'$ and $S'\cap S=\emptyset$, and $\Aut(\Sigma)_1$ is transitive on $S'$. 

We claim that also $S$ is an orbit of $\Aut(\Sigma)_1$. 
The set $S\cup\{1\}=X\cup Y$ is invariant under $\Aut(\Sigma)_1$, and in the proof of Step 3 we noted that the subgroup $\Aut(X)\times \Aut(Y)$ of $\Aut(\Sigma)_1$ is transitive on each of $X\setminus\{1\}$ and $Y\setminus\{1\}$. Moreover we are assuming that $\Aut(\Sigma)$ is transitive on $V(\Sigma)$ and hence, since $\Sigma$ is locally $2$-arc-transitive, $\Aut(\Sigma)$ is transitive on the arcs of $\Sigma$. Thus $\Aut(\Sigma)$ contains an element $\sigma$ which maps the arc $(X,Y)$ of $\Sigma$ to the arc $(Y,X)$. Since $\Sigma$ is the clique graph of $\G$ (Lemma~\ref{lem:prop-mixed-dih}(1)), $X, Y$ (as subsets of $\newH(n)$) are maximal cliques of $\G$ and are interchanged by $\sigma$. In particular, $\sigma$ induces an automorphism of the subgraph of $\Gamma$ induced on $X\cup Y$. The identity $1$ is adjacent in $\G$ to every vertex of $S$, while each other vertex $z\in S$ is adjacent to only $|X|-1$ elements of $S\cup\{1\}$. Thus $\sigma$ must fix $1$ and interchange $X\setminus\{1\}$ and $Y\setminus\{1\}$. It follows that  $\Aut(\Sigma)_1$ is transitive on $S$, proving the claim. Thus $\Aut(\Sigma)$ has exactly two orbits on the arcs of $\Lambda$, namely the arcs $(w,z)$ with $wz^{-1}\in S$ and those with $wz^{-1}\in S'$.

Next we identify certain small subgraphs of $\Lambda$. For any $[x,y]\in S'$ and $y'\in Y\setminus\{1\}$, by Theorem~\ref{prop-2groupNew}~(4) we have $[[y,x],y']=1$ and $[[x,y],y']=1$, so $[x,y]y'=y'[x,y]$, and
\[
(1, y', [x,y]y', [x,y], 1)
\]
is a $4$-arc of $\Lambda$. Moreover, since $[x,y]y'\not\in S\cup S'$, it follows that the subgraph of $\Lambda$ induced on $C(x,y,y'):=\{1, y', [x,y]y', [x,y]\}$ is a $4$-cycle. 

Now we choose $x'=x_2\in X_0$ and $a=x_1\in X_0$, $b=y_1\in Y_0$ so that $[a,b]\in S'$. These elements arise as images under $\sigma$ as follows: there exists $[x,y]\in S'$ such that $[x,y]^\s=[a, b]$, and there exists $y'\in Y$ such that $(y')^\s=x'$. Thus the subgraph of  $\Lambda$ induced on $C':= C(x,y,y')^\sigma = \{ 1, x', ([x,y]y')^\s, [a,b]\}$ is a $4$-cycle including the $2$-arc $([a,b],1,x')$. Thus, setting $z:= ([x,y]y')^\s$, this $4$-cycle is $(1, x', z, [a,b], 1)$ and hence $z= sx'=t[a,b]$ for some $s,t\in S\cup S'$. As these four vertices are pairwise distinct, $sx'\neq1$ and $tsx'=[a,b]\neq1$.

Each of $t,s$ lies in either $S$ or $S'$, giving four possible combinations. We obtain a contradiction from each possibility as follows. First, if both $t,s\in S'$, then $x'=st[a,b]$ and $st[a,b]\in\newH(n)'$ while $x'\in X\setminus\{1\}$, and we have a contradiction since by \cite[Lemma 2.1]{HPZ}, $\phi(x')\ne1$ while $\phi(st[a,b])=1$.   Next suppose that $t\in S'$ and $s\in S$. Then $sx'=t[a,b]\in\newH(n)'$ and hence $\phi(sx')=\phi(t[a,b])=1$, which implies that $sx'=1$ (by \cite[Lemma 2.1]{HPZ}),  a contradiction. Thirdly, suppose that $s,t\in S$. Then $tsx'=[a,b]\in\newH(n)'$ and hence $\phi(tsx')=\phi([a,b])=1$. Again we conclude that $tsx'=1$ by \cite[Lemma 2.1]{HPZ}, which is a contradiction. 

This leaves the case $t\in S, s\in S'$, and hence $\phi(x')=\phi(sx')=\phi(t[a,b])=\phi(t)$. Then by \cite[Lemma 2.1]{HPZ}, we must have $t\in X\setminus\{1\}$ and $\phi(x')=\phi(t)$ implies that $t=x'$. Thus  $x'sx'=[a,b]$ and so $s=x'[a,b]x'=[a,b]^{x'}\in S'$, which implies that $[a,b][[a,b],x']=[a,b]^{x'}=s=[c,d]$ for some $c\in X\setminus\{1\}$ and $d\in Y\setminus\{1\}$. 
Now $c=x_{i_1}x_{i_2}\dots x_{i_k}$ and $d=y_{j_1}y_{j_2}\dots y_{j_\ell}$ for some $x_{i_1},x_{i_2},\dots, x_{i_k}\in X_0$ and $y_{j_1}, y_{j_2},\dots, y_{j_\ell}\in Y_0$, and we choose these expressions with $k, \ell$ minimal. We now apply Theorem~\ref{prop-2groupNew}~(6) to $[c,d]$, observing that $[a,b]\newH(n)_3= [c,d]\newH(n)_3$, recalling that $a=x_1\in X_0$ and $b=y_1\in Y_0$), and using the fact that the cosets $[x_{i_u},y_{j_v}]\newH(n)_3$, for $1\leq u,v\leq n$, form a basis for $\newH(n)'/\newH(n)_3\cong C_2^{n^2}$ (see Theorem~\ref{prop-2groupNew}~(3)). We deduce that $k=\ell=1$, so that $a=c=x_{i_1}=x_1$ and $b=d=y_{j_1}=y_1$. The equality $[a,b][[a,b],x']=[c,d]$ then implies that $[[a,b],x']=1$, that is to say, $[[x_1,y_1],x_2]=1$, which contradicts Theorem~\ref{prop-2groupNew}~(3). Thus $\Aut{\Gamma}$ acts intransitively on the vertices of $\Gamma$, and $\Gamma$ is semisymmetric, completing the proof of Theorem~\ref{t:semi-sym}.
\hfill\qed

\begin{rem}\label{rem-magma}
{\rm To prove $\Sigma$ is semisymmertic in case $n=2$, we make use of a Magma~\cite{BCP} computation, which we now describe. First, the group $\newH(2)$ is input in the category GrpFP via the presentation given in Definition~\ref{def:Hn}. Next, the {\tt pQuotient} command is used to construct the largest $2$-quotient $H2$ of $\newH(2)$ having lower exponent-$2$ class at most 100 as group in the category GrpPC. Comparing the orders of these groups, we find $|\newH(2)|=|H2|$, so that $\newH(2)\cong H2$. Next we construct the graph $\Sigma$. Computation shows that $\Sigma$ is edge-transitive but not vertex-transitive, and has valency $4$. We have made available the Magma programs in the following Appendix.}
\end{rem}

\section*{Appendix: Magma programs used in the proof of Theorem~\ref{t:semi-sym} in the case $n=2$.}

\f{Input the group $\newH(2)$:}\smallskip
\\ 
{\tt G<x1,x2,y1,y2>:=Group<x1,x2,y1,y2|
x1\textasciicircum2, x2\textasciicircum2, y1\textasciicircum2, y2\textasciicircum2,
\\
(x1,x2)=(y1,y2)=1,
\\
(x1,y1)\textasciicircum2=(x1,y2)\textasciicircum2=(x2,y1)\textasciicircum2=(x2,y2)\textasciicircum2=1,
\\
((x1,y1),x2)\textasciicircum2=((x1,y1),y2)=1,
\\
((x1,y2),x2)\textasciicircum2=((x1,y2),y1)=1,
\\
((x2,y1),x1)\textasciicircum2=((x2,y1),y2)=1,
\\
((x2,y2),x1)\textasciicircum2=((x2,y2),y1)=1,
\\
(x1,((x1,y1),x2))=(x2,((x1,y1),x2))=(y1,((x1,y1),x2))=(y2,((x1,y1),x2))=1,
\\
(x1,((x1,y1),y2))=(x2,((x1,y1),y2))=(y1,((x1,y1),y2))=(y2,((x1,y1),y2))=1,
\\
(x1,((x1,y2),x2))=(x2,((x1,y2),x2))=(y1,((x1,y2),x2))=(y2,((x1,y2),x2))=1,
\\
(x1,((x1,y2),y1))=(x2,((x1,y2),y1))=(y1,((x1,y2),y1))=(y2,((x1,y2),y1))=1,
\\
(x1,((x2,y1),x1))=(x2,((x2,y1),x1))=(y1,((x2,y1),x1))=(y2,((x2,y1),x1))=1,
\\
(x1,((x2,y1),y2))=(x2,((x2,y1),y2))=(y1,((x2,y1),y2))=(y2,((x2,y1),y2))=1,
\\
(x1,((x2,y2),x1))=(x2,((x2,y2),x1))=(y1,((x2,y2),x1))=(y2,((x2,y2),x1))=1,
\\
(x1,((x2,y2),y1))=(x2,((x2,y2),y1))=(y1,((x2,y2),y1))=(y2,((x2,y2),y1))=1>;}
\smallskip

\f {Construct the largest 2-quotient group of $\newH(2)$ having lower exponent-$2$ class at most 100 as group in the category GrpPC:}
\smallskip
\\
{\tt H2,q:=pQuotient(G,2,100);}\smallskip

\f{Order of $\newH2$ (The result shows that $|H2|=|\newH(2)|$, and so $H2\cong \newH(2)$):}\smallskip
\\
{\tt FactoredOrder(H2);}\smallskip

\f{Construct the graph $\Sigma$:}
\smallskip
\\
{\tt X:=sub<H2|x1,x2>;
Y:=sub<H2|y1,y2>;
\\
Vsigma1:=$\{\}$;
\\
for g in H2 do
\\
Xg:=$\{\}$;
\\
for a in X do
\\
Include($^\sim$Xg, a*g);
\\
end for;
\\
Include($^\sim$Vsigma1,Xg);
\\
end for;
\\
Vsigma2:=$\{\}$;
\\
for g in H2 do
\\
Yg:=$\{\}$;
\\
for b in Y do
\\
Include($^\sim$Yg, b*g);
\\
end for;
\\
Include($^\sim$Vsigma2,Yg);
\\
end for;
\\
Vsigma:=Vsigma1 join Vsigma2;
\\
Esigma:=$\{\{$x,y$\}$: x in Vsigma1, y in  Vsigma2 | $\sharp$(x meet y) ne 0$\}$;
\\
Sigma:=Graph<Vsigma|Esigma>;
}

\smallskip
\f {Test if $\Sigma$ is a tetravalent semisymmetric graph:}
\smallskip
\\
{\tt 
IsVertexTransitive(Sigma);
\\
IsEdgeTransitive(Sigma);
\\
Valence(Sigma);
}

\section*{Acknowledgements}

The first author has been supported by the Croatian
Science Foundation under the project 6732. The second author is grateful for Australian Research Council Discovery Project Grant DP230101268. The third author was supported by the National Natural Science Foundation of China (12071023, 12161141005) and the 111 Project of China (B16002).


\begin{thebibliography}{99}

\bibitem{BCP}
W. Bosma, J. Cannon, C. Playoust,
The MAGMA algebra system I: The user language,
{\em J. Symbolic Comput.} 24 (1997) 235--265. \\[-20pt]


\bibitem{C-M-M-P}
M.~Conder, A.~Malni\v{c}, 
D.~Maru\v{s}i\v{c}, and P.~Poto\v{c}nik,
A census of semisymmetric cubic graphs on up
to $768$ vertices, {\em J Algebr. Comb.} {\bf 23} (2006), 255--294.

\bibitem{C-Z-F-Z}
M. Conder, J.-X. Zhou, Y.-Q. Feng, M.-M. Zhang, Edge-transitive bi-Cayley graphs, {\em J. Comb. Theory B}  145 (2020) 264--306.\\[-20pt]





\bibitem{FanLLP2013}
W. Fan, D. Leemans, C.H. Li, J. Pan, Locally $2$-arc-transitive complete bipartite graphs, {\em J. Comb.  Theory A} 120 (2013) 683--699. \\[-20pt]

\bibitem{Fo}
J.~Folkman, 
Regular line-symmetric graphs.
{\em J. Combinatorial Theory} {\bf3} (1967), 215--232.


\bibitem{GAP4}%
The GAP~Group.
\newblock {\em {GAP -- Groups, Algorithms, and Programming, Version 4.11.1}}, 2021.\\[-20pt]


\bibitem{GLP-tamc}
M. Giudici, C.H. Li, C.E. Praeger, Analysing finite locally $s$-arc transitive graphs, {\em Tran. Amer. Math. Soc.} 356 (2004) 291--317.\\[-20pt]

\bibitem{Godsil-1981}
C.D. Godsil, On the full automorphism group of a graph, {\em Combinatorica} 1 (1981) 243--256. \\[-20pt]



\bibitem{Hall-book}
M. Hall, The Theory of Groups, The Macmillan Company, New York.


\bibitem{HPZ}
D.R. Hawtin, C.E. Praeger, J.-X. Zhou, A characterisation of edge-affine $2$-arc-transitive covers of $\K_{2^n,2^n}$. Preprint, arXiv:2211.16809 (2022).\\[-20pt]


\bibitem{Huppert}
B. Huppert, Eudiche Gruppen I, Springer-Verlag, New York, 1967.

\bibitem{IP}
A.A. Ivanov, C.E. Praeger, On finite affine $2$-arc-transitive graphs, {\em Eur. J. Combin.} 14 (1993) 421--444.

\bibitem{AVIv}
A.V. Ivanov, 
On edge but not vertex transitive regular graphs. {\em Combinatorial design theory}, pp. 273--285,
North-Holland Math. Stud., 149, Ann. Discrete Math., 34, North-Holland, Amsterdam, 1987.




\bibitem{L-BLMS-2001}
C.H. Li, Finite $s$-arc transitive graphs of prime-power order, {\em Bull. London Math. Soc.} 33 (2001) 129--137. \\[-20pt]

\bibitem{L-TAMS-2006}
C.H. Li, Finite edge-transitive Cayley graphs and rotary Cayley maps, {\em Tran. Amer. Math. Soc.} 358 (2006) 4605--4635. \\[-20pt]

\bibitem{LMP2009}
C.H. Li, L. Ma, J. Pan, Locally primitive graphs of prime-power order, {\em J. Aust. Math. Soc.} 86 (2009) 111--122. \\[-20pt]

\bibitem{Li-Pan}
C.H. Li, J.M. Pan, Finite $2$-arc-transitive abelian Cayley graphs, {\em Eur. J. Comb.} 29 (2008) 148--158.\\[-20pt]





\bibitem{Imp}
C.E. Praeger, Imprimitive symmetric graphs, {\em Ars Combinatoria} 19A (1985) 149--163.

\bibitem{LPVZ}
C.H. Li, C.E. Praeger, A. Venkatesh and S. Zhou, Finite locally-quasiprimitive graphs, {\em Discrete Math.} 246 (2002), 197--218.

\bibitem{ONS}
C.E. Praeger, An O'Nan--Scott Theorem for finite quasiprimitive permutation groups, and an application to 2-arc transitive graphs, {\em J. London Math. Soc.}(2) 47 (1992), 227--239.

\bibitem{Bip}
C.E. Praeger, On a reduction theorem for finite bipartite $2$-arc transitive graphs, {\em Australas. J.
Combin.} 7 (1993) 21--36.


\bibitem{PS}
C.E. Praeger, C. Schneider, {\em Permutation Groups and Cartesian Decompositions}, LMS Lecture Series 449, Cambridge University Press, Cambridge, 2018.

\bibitem{Robinson}
D.J. Robinson, A Course in the Theory of Groups, Second Edition, Springer, New York, 1996.\\[-20pt]

\bibitem{T47}
W.T. Tutte, A family of cubical graphs, {\em Proc. Cambridge Philos. Soc.} 43 (1947), 459--474.



\bibitem{ZF}
J.-X.~Zhou, Y.-Q.~Feng,
The automorphisms of bi-Cayley graphs,
{\em J. Comb. Theory B} 116 (2016) 504--532.\\[-20pt]

\end{thebibliography}
\end{document}